\documentclass[10pt]{amsart}

\usepackage[latin1]{inputenc}
\usepackage{amsmath,amsthm} 
\usepackage{amssymb, latexsym}   

\newtheorem{theorem}{Theorem}

\newtheorem{lemma}{Lemma}

\begin{document}

\title[Inverse Lévy subordinator]{Scaling limit results for the sum of
  many inverse Lévy subordinators}

\author{Ingemar Kaj \and  Anders Martin-Löf}
\address{IK, Dept.\ of Mathematics, Uppsala
University, Box 480, SE 751 06 Uppsala, Sweden}
\email{ikaj@math.uu.se}

\address{AM-L, Mathematical Statistics, Stockholm University, SE 106
91 Stockholm, Sweden}
\email{andersml@math.su.se}

\date{Revised version, March 2012}

\begin{abstract}

The first passage time process of a Lévy subordinator with
heavy-tailed Lévy measure has long-range dependent paths. The random
fluctuations that appear under two natural schemes of summation and
time scaling of such stochastic processes are shown to converge
weakly. The limit process is fractional Brownian motion in one case
and a non-Gaussian and non-stable process in the other case.  The
latter appears to be of independent interest as a random process that
arises under the influence of coexisting Gaussian and stable domains
of attraction and is known from other applications to provide a bridge
between fractional Brownian motion and stable L\'evy motion.
\end{abstract}
\keywords{Long-range dependence, weak convergence, fractional
Brownian motion.}

\maketitle

\section{Introduction and statement of results}

A Lévy subordinator $\{X_t, t\ge 0\}$ is a real-valued random process
with independent and stationary increments and increasing pure-jump
trajectories. The inverse process $\{T_x,x\ge 0\}$ defined by the
first passage times $T_x=\inf\{t\ge 0: X_t>x\}$ has nondecreasing,
 trajectories, where the lengths of the flat pieces of $\{T_x\}$
correspond to the jump sizes of $\{X_t\}$. The dependence structure in
the paths of the inverse process is entirely different from that of
the Lévy subordinator, since big jumps in the Lévy process may cause
strong dependencies that last over a considerable period of evolution
of the path of its inverse. In this paper we take a scaling approach
to study the nature of the random fluctuations that build up as a
result of such long-memory effects. By superposing a large number of
paths of the inverse Lévy process and simultaneously scale the time
parameter of the process, we obtain scaling limit results for the
centered and normalized superposition process.

In somewhat more detail, our starting point is a Lévy subordinator
with Lévy measure $\nu(dx)$ of regularly varying tail with index
$1+\beta$, $0<\beta<1$. In particular, $\mu:=\int x\nu(dx)<\infty$.
The initial distribution of the subordinator process is chosen such
that the resulting inverse process has stationary increments and
expected value $E(T_x)=x/\mu$.  Letting $\{T_x^i\}_{i\ge 1}$ be a
collection of independent copies of $\{T_x\}$, our main result is the
derivation of a limit process for the summation scheme
\[
\frac{1}{a}\sum_{i=1}^m (T_{ax}^i -\frac{1}{\mu}ax),\quad x\ge 0,
\]  
as both $m$ and $a=a_m$ tend to infinity in such a way that $m$ is of
the same order of magnitude as $a^\beta$, modulo slowly varying
functions. The reason for this choice of scaling is to attempt to
trace the superposition process on a time scale that captures the size
of the fluctuations around its mean.  In the asymptotic limit appears
a non-Gaussian, non-stable process with long-range dependence, which
is known to arise also in other related models and has been called
fractional Poisson motion, \cite{gaigalaskaj}, \cite{gaigalas},
\cite{kajtaqqu}, \cite{dombrykaj2011}. The general study
\cite{dombrykaj2012} of higher-order
moment measures for heavy-tailed renewal point processes, provides a
unified framework of \cite{gaigalas} and the present work. 

To give a heuristic context for the topics of interest in this work,
let us recall the following limit result for Lévy processes. Writing
$\alpha=1+\beta$, the centered and scaled process $(X_t-\mu
t)/t^{1/\alpha}$ converges in distribution as $t\to\infty$ to a random
variable $Z_\alpha$, having a stable distribution with stable index
$\alpha$. If we write $\Gamma_x$ for the overshoot at $x$, so that
$X_{T_x}=x+\Gamma_x$, then
\[
{T_x-x/\mu\over x^{1/\alpha}}=- {X_{T_x}-\mu T_x\over T_x^{1/\alpha}} 
        \Big({T_x\over x}\Big)^{1/\alpha}{1\over\mu}
       +{\Gamma_x\over\mu x^{1/\alpha}}.
\]
In this relation, $T_x/x\to 1/\mu$ as $x\to\infty$ by the law of large
numbers. It can be shown moreover that the second term on the right
hand side is a remainder term with $\Gamma_x/x^{1/\alpha}\to 0$ as
$x\to\infty$. Therefore $(T_x-x/\mu)/x^{1/\alpha}$ converges in
distribution to $-Z_\alpha/\mu^{1+1/\alpha}$ as $x\to\infty$.
Proceeding heuristically, with $m\sim a^\beta$ we may rewrite
the superposition process either as
\[
{1\over a}\sum_{i=1}^m (T_{ax}^i -{1\over\mu}ax)
 \sim   {1\over a^{1-\beta/2}}\int_0^{ax}
{1\over m^{1/2}}\sum_{i=1}^m (dT_u^i-{1\over\mu}du)
\]  
or
\[
{1\over a}\sum_{i=1}^m (T_{ax}^i -{1\over\mu}ax)
\sim {1\over m^{1/(1+\beta)}}\sum_{i=1}^m 
        {T_{ax}^i -ax/\mu\over a^{1/(1+\beta)}}. 
\]
The first representation emphasizes a sequence of random variables in
the domain of attraction of a Gaussian law ($m\to\infty$ with $a$
fixed). The second representation highlights a sequence in the domain
of attraction of a stable law with index $1+\beta$ ($a\to\infty$ with
$m$ fixed), which is the type of convergence just discussed above. For
the limit regime of interest in our case Gaussian and stable
attraction appear to coexist and both influence the resulting limit
process.

The main result (Theorem 2 below) is a scaling limit theorem for the
intermediate type rescaling regime indicated above. In parallel to
this we discuss the scaling regime of Gaussian predominance, leading
to fractional Brownian motion in the limit (Theorem 1). Scaling limit
results with fractional Brownian fluctuations are known for a variety
of models, such as modeling random variation in aggregated data
traffic streams. For an introduction and overview of these topics and
discussion of the modeling context, as well as detailed statements and
derivations of such results, see e.g.\ \cite{taqqu, willingeretal,
  kajtaqqu}.

The model is introduced in detail and all results are stated in
Section 1 of the paper.  We then focus on the proof of Theorem 2 for
the intermediate scaling regime, starting with the analysis of
marginal distributions in Section 2. The main technique we use for the
study of the one-dimensional distributions of the scaled processes and
their limit behavior is that of double transforms in the sense of
taking Laplace transforms in the time variable of the logarithmic
moment generating function of the random variables. In Section 3 we
continue with a study of the finite-dimensional distributions, which are
obtained from recursive sets of integral equations for the
finite-dimensional cumulant functions. Finally in Section 4 we provide
a summary of the proof of Theorem 1 for Gaussian scaling, where each
parallel step turns out to be simpler, and the proof of tightness. 

\subsection{A Lévy subordinator and its inverse}

Let $\{\widetilde X_t,t\ge 0\}$, $\widetilde X_0=0$, denote a Lévy
subordinator with right-continuous paths, having drift zero and Lévy
measure $\nu(a,b)=\int_a^b\,\nu(dx)$ with no atom at zero, such that
\begin{equation}\label{muassumption}
\int_0^\infty (1\wedge x)\nu(dx)<\infty 
\quad \mbox{and}\quad 
\mu=\int_0^\infty x\nu(dx)<\infty,
\end{equation}
which implies that the first moment is finite, $E(\widetilde X_t)=\mu
t<\infty$.  The Laplace transform is
given by $-\ln E(e^{-u\widetilde X_t})=t \Phi(u)$, $u\ge 0$,
with Laplace exponent
\[
\Phi(u)=\int_0^\infty(1-e^{-ux})\,\nu(dx).
\]
Let $X_t=X_0+\widetilde X_t$ denote the corresponding delayed
subordinator process with general initial distribution $X_0$ assumed
to be independent of $\{\widetilde X_t\}$.  We will study the case
when $X_0>0$ has distribution function
\begin{equation}\label{X_init}
P(X_0\le x)={1\over \mu}\int_0^x \int_y^\infty \nu(ds)\,dy,
\end{equation}
for which $E(e^{-u X_0})
={1\over \mu u}\Phi(u)$
and so
\begin{equation}\label{mgf_X}
E(e^{-u X_t})={1\over \mu u}\Phi(u)
\,\exp\{-t\Phi(u)\}\quad u\ge 0.
\end{equation}

Next we introduce the first passage process of the
subordinator. Useful references are Bertoin \cite{bertoin},
\cite{bertoin2}. Van Harn and Steutel, \cite{vanharnsteutel},
investigate stationarity properties of delayed subordinators and
derive closely related results to those in Lemma 1 and Lemma 3 below.
The entrance time of the Lévy process $\{X_t\}$ into a set $B$ is
defined by $T_B=\inf\{t\ge 0: X_t\in B\}$. For any open set $B$, $T_B$
is a stopping time. The first passage time $T_x=T_{(x,\infty)}$
strictly above a level $x$ is the entrance time into $(x,\infty)$,
that is
\[
T_x=\inf\{t\ge 0: X_t>x\}, \quad x\ge 0,
\]
which is a right-continuous function with left limits.  Since
$X_t\uparrow\infty$ as $t\uparrow\infty$, we have $T_x<\infty$ for all $x$
and $P(T_x\le t)=P(X_t>x)$, $x,t\ge 0$. 
Also, 
\[
E(T_x)=\int_0^\infty P(T_x>t)\,dt
  =E\int_0^\infty 1_{\{X_t\le x\}}\,dt.
\]
Hence 
\[
\int_0^\infty ue^{-ux} E(T_x)\,dx=\int_0^\infty E(e^{-uX_t})\,dt
  = {1\over \mu u}
\]
in view of (\ref{mgf_X}), and therefore
\[
E(T_x)={1\over \mu}x.
\]
We call $\{T_x\}$ the inverse Lévy subordinator and the process
\[
\widetilde T_x=\inf\{t\ge 0: \widetilde X_t>x\}, \quad x\ge 0,
\]
the pure inverse Lévy subordinator. 
The path-regularity of
$\{T_x\}$ is determined by the distribution of the small jumps of
$\{X_t\}$, manifest in the asymptotic behavior of $\nu(dx)$ for $x$
close to 0. The nature of the paths varies considerably with compound
Poisson processes as one extreme case.  These are the subordinators
for which the Lévy measure is finite on the positive half line and the
passage time process has piece-wise constant trajectories of lengths
drawn from the probability distribution $\nu(0,x]/\nu(0,\infty)$,
$x\ge 0$, and with exponentially distributed jumps.  On the other
hand, if the number $\sigma=\sup\{\alpha>0:\lim_{\lambda\to\infty}
\lambda^{-\alpha}\Phi(\lambda)=\infty\}$, known as the lower index of
the subordinator, is positive (and $\le 1$ because of assumption
(\ref{muassumption})) then the inverse process $\{T_x\}$ is a.s.\
$\gamma$-Hölder continuous on any compact interval, for each index
$\gamma<\sigma$, see \cite{bertoin}, Ch.\ 3 ($X_0=0$). 

The scaling problem studied in this work involves weak convergence in
the sense of convergence of finite-dimensional distributions plus a
tightness property in the space $C=C[0,\infty)$ of continuous random
  processes.  The path space is such that each $C[0,T]$, $T>0$, is
  equipped with the topology of convergence in supremum norm.

We will prove below the following 
\begin{lemma}\label{statincr}
The inverse subordinator process $\{T_x, x\ge 0\}$ has stationary
increments.
\end{lemma}

\subsection{Scaling limit theorem}

Our basic assumption is that the Lévy measure $\nu$ is regularly
varying at infinity with index $1+\beta$, $0<\beta<1$, i.e.
\begin{equation}\label{assumpbeta}
\int_x^\infty \nu(dy)\sim {1\over x^{1+\beta}}L(x),\quad x\to\infty,
\end{equation}
where $L$ is a slowly varying function and we write $f(x)\sim g(x)$ if
$f$ and $g$ are positive functions and $f(x)/g(x)\to 1$ as
$x\to\infty$.  The summation schemes to be applied involve speeding up
the time parameter using a rescaling sequence $a_m\to\infty$, either
such that
\begin{equation}\label{FCR}
m L(a_m)/a_m^\beta\to \infty,\quad m\to\infty,
\end{equation}
or such that
\begin{equation}\label{ICR}
m L(a_m)/a_m^\beta\to c^\beta\,\mu,\quad m\to\infty,
\end{equation}
where $c$, $0<c<\infty$, is an additional parameter that signifies the
relative change of scales of size and time. In addition, we assume
that the lower index of the Lévy measure $\nu$ satisfies 
\begin{equation}\label{lowerindex}
\sigma=\sup\{\alpha>0:\lim_{\lambda\to\infty}
\lambda^{-\alpha}\Phi(\lambda)=\infty\}>\beta.
\end{equation}

\begin{theorem} \label{thmfbm}
Under assumptions (\ref{assumpbeta}) and (\ref{lowerindex}), let $a_m$ be a sequence such
that $a_m\to\infty$ as $m\to\infty$ and (\ref{FCR}) holds, and define
$b_m$ by 
\begin{equation}\label{FCR_b}
b_m^2=m a_m^{2-\beta}L(a_m)/\mu. 
\end{equation}
Then, in the sense of weak convergence of random processes in $C$,
\begin{equation}\label{FCR_conv}
\Big\{{1\over b_m}\sum_{i=1}^m (T^i_{a_m x}-{1\over\mu} a_m x),\;x\ge 0\Big\}
\Rightarrow   \{\mu^{-1}\sigma_\beta B_H(x),\; x\ge 0\},
\end{equation}
where 
\[
\sigma^2_\beta={2\over\beta(1-\beta)(2-\beta)},\quad H=1-\beta/2,
\]
and $B_H$ is standard fractional Brownian motion with Hurst index $H$,
i.e.\ the Gaussian process with stationary increments, variance
$V(B_H(t))=t^{2H}$ and continuous sample paths. 

\begin{eqnarray}\nonumber
\log E \exp \Big\{\sum_{i=1}^n \theta_i B_H(x_i)\Big\}
={1\over 4}\sum_{i=1}^n \sum_{j=1}^n \theta_i \theta_j
(x_i^{2-\beta}+x_j^{2-\beta}-(x_i-x_j)^{2-\beta}), 
\label{logmgf_FBM}
\end{eqnarray}
where $0=x_0\le x_1\le \dots \le x_n$, $n\ge 1$. 
\end{theorem}

\begin{theorem} \label{thmICR}
Under the assumptions (\ref{assumpbeta}) and (\ref{lowerindex}), if $a_m$ is a sequence such
that $a_m\to\infty$ and (\ref{ICR}) holds for some constant $c>0$ as
$m\to\infty$, then 
\begin{equation}\label{ICRconv}
\Big\{{1\over a_m}\sum_{i=1}^m (T^i_{a_m x}-{1\over\mu} a_m x),\;x\ge 0\Big\}
\Rightarrow   \{-\mu^{-1}c\,Y_\beta(x/c),\; x\ge 0\},
\end{equation}
in the sense of weak convergence in $C$. Here $\{Y_\beta(x),x\ge 0\}$
is a zero mean stochastic process with continuous paths and
finite-dimensional distributions characterized by the cumulant
generating function
\begin{eqnarray} \label{mgfndim}
\lefteqn{\hspace{-10mm}
\log E \exp \Big\{\sum_{i=1}^n \theta_i
\big(Y_\beta(x_i)-Y_{\beta}(x_{i-1}) \big) \Big\}
=\frac{1}{\beta} \sum_{i=1}^n
\theta_i^2 \int_0^{\Delta x_i}\!\int_0^v e^{\theta_i u}
u^{-\beta}\,dudv \nonumber 
}\\ 
&&\hspace{-20pt}\quad 
+ \frac{1}{\beta}
\sum_{i=1}^{n-1} \sum_{j=i+1}^{n} \theta_i \theta_j
\exp\Big\{\sum_{k=i+1}^{j-1} \theta_k \Delta x_k\Big\} 
\nonumber \\ 
&& \qquad\times 
\int_0^{\Delta x_i}\! 
\int_0^{\Delta x_j}e^{\theta_j u} e^{\theta_i v}
(x_{j-1}-x_i+u+v)^{-\beta}\,dudv, 
\end{eqnarray}
where $0=x_0\le x_1\le \dots \le x_n$, and $\Delta x_i=x_i-x_{i-1}$,
$i=1,\dots,n$.
\end{theorem}

\paragraph{{\bf Remarks}} 

{\bf a)} It is known that $\{Y_\beta\}$ has a simple representation as the
stochastic integral
\[
Y_\beta(x)=\int_{R\times [0,\infty)}\int_0^x 1_{\{s<y<s+v\}}\,dy \, \widetilde
N(ds,dv)
\]
where $\widetilde N(ds,dv)=N(ds,dv)-n(ds,dv)$ and $N(ds,dv)$ is a
Poisson random measure on $R\times [0,\infty)$ with intensity measure 
$n(ds,dv)=(1+\beta)ds\, v^{-2-\beta}dv$, \cite{kaj, gaigalas, kajtaqqu}.  
We have not been able to find however a method of proof of the present
results which utilizes more directly this inherent Poisson structure
of the model. 

{\bf b)} The process $\{Y_\beta(x)\}$ has been derived in 
\cite{gaigalaskaj} as a limit process in the setting of a
superposition of independent renewal processes with stationary
increments and heavy-tailed inter-renewal distribution, and in
\cite{kaj} and 
\cite{kajtaqqu} for an infinite source Poisson process with heavy-tailed
activity periods. The motivation is partly from modeling
the total traffic load generated by many independent sources at an
arrival point in a data traffic network.  In these references
condition (\ref{FCR}) is called {\it fast connection rate} and
(\ref{ICR}) {\it intermediate connection rate}. They are compared to
an alternative third scaling regime of {\it slow} connection rate, for
which the limit process turns out to be a stable Lévy process with
stable index $\alpha=1+\beta$, see also {\it et al.}
\cite{mikoschetal} or {\it et al.} \cite{willingeretal}.

\smallskip
{\bf c)} Proofs of the following properties among others can be found
in Gaiga\-las and Kaj \cite{gaigalaskaj}. The process $\{Y_\beta\}$
has stationary increments and continuous trajectories.  The process is
not self-similar. The higher moments are of the order
$E(Y^k_\beta(x))\sim \mbox{\rm const}\, x^{k-\beta}$, $k\ge 2$, for
large $x$. Specifically, the second-order properties (mean, variance,
covariance) are the same (modulo constants) as those for fractional
Brownian motion, where\-as higher order moments are different. For
example, $\{Y_\beta\}$ is positively skewed. The paths are
$\gamma$-Hölder continuous for all $\gamma<1-\beta/2$ (not $\gamma<1$
as claimed in \cite{gaigalaskaj}).


\smallskip
{\bf d)} 
The renewal processes studied in \cite{gaigalaskaj} can
be viewed as discrete local time processes of discrete regenerative
sets (ranges of compound Poisson subordinators). In this light, the
present situation is the natural analogue for continuous local time
processes of perfect regenerative sets (ranges of subordinators that
are not compound Poisson). One can expect the scaling limits to
transfer since they are large-time asymptotics which should not depend
on the local structure. Some relevant references for the connections
of regenerative sets and subordinators are \cite{fristedt},
and \cite{gnedinpitman}.

\section{Analysis of the marginal distrbution}

As a preliminary for the proof of Theorem 2 we observe the following
properties of the functions introduced in (\ref{mgfndim}), which are
straightforward to verify.
\begin{lemma}\label{lemmaconsistency}
Relation (\ref{mgfndim}) defines a consistent family of
finite-dimensional distributions, such that for any $c>0$
\begin{eqnarray*}
&&\log E \exp \Big\{\sum_{i=1}^n \theta_i
(cY_\beta(x_i/c)-cY_{\beta}(x_{i-1}/c)) \Big\}\\
  &&\quad =c^\beta\log E \exp \Big\{\sum_{i=1}^n \theta_i
 (Y_\beta(x_i)-Y_{\beta}(x_{i-1}))\Big\}.
\end{eqnarray*}
\end{lemma}

The main part of the proofs of Theorem 1 and Theorem 2 consists in
establishing convergence of the scaled $n$-point cumulant functions
\begin{eqnarray}\nonumber
\lefteqn{\log E\exp\Big\{ \sum_{i=1}^n \theta_i  
 {1\over b_m}\sum_{k=1}^m (T^{(k)}_{a_m x_i}-{1\over\mu}a_m x_i)
       \Big\}}\\ 
&& = m E\Big[\exp\Big\{\sum_{i=1}^n
  {\theta_i \over b_m}(T_{a_m x_i}-{1\over\mu}a_m x_i)\Big\}-1\Big]
    +{\mathcal O}(1/m)  \label{defmgf}
\end{eqnarray}
toward the corresponding functionals of the limit processes. As a
preparation we study the joint distribution $(T_x,\Gamma_x)$, where
$\{\Gamma_x\}$ is the overshoot process, and other properties of the
one-dimensional marginal distributions of $T_x$.

\subsection{Marginal distributions and the overshoot process}

The overshoot process $\{\Gamma_x, x\ge 0\}$ associated with the first
passage time $\{T_x\}$ is defined for $x\ge 0$ by 
\[
\Gamma_x=X_{T_x}-x,
\]
and represents at time $x$ the remaining time until the next point of
increase of the inverse subordinator.  The following Lemma is a
special case of a result valid for general Lévy processes adaptated to
the case of a general initial distribution $X_0$. For a proof see
Theorem 49.2 in Sato \cite{sato}. 
\begin{lemma} \label{jointgammaT} 
For $u>0$, $\theta<\Phi(u)$, and $v>0$ with $v\not= u$,
\[
\int_0^\infty ue^{-ux} E(e^{\theta T_x-v \Gamma_x})\,dx
  = {u\over u-v}\left({\Phi(v)\over\mu v} -
        {\Phi(v)-\theta\over \Phi(u)-\theta} \,{\Phi(u)\over\mu u}
\right).
\]
\end{lemma}

\paragraph{Proof of Lemma \ref{statincr}}

In Lemma \ref{jointgammaT}, take $u>0$ and $v>0$, $u\not= v$, and let
$\theta=0$. We obtain 
\[
\int_0^\infty ue^{-ux} E(e^{-v \Gamma_x})\,dx
  = {u\over u-v}\left({\Phi(v)\over\mu v} -
        {\Phi(v)\over \Phi(u)} \,{\Phi(u)\over\mu u} \right)
  = {\Phi(v)\over\mu v}.
\] 
Hence, for any $x\ge 0$, $\Gamma_x \stackrel{d}{=} X_0$. Consequently,
for each $x$ the increment process $T_{x+y}-T_x$, $y\ge 0$, begins
with a flat period for a duration of time having the distribution
$X_0$, which is just the same behavior as the original process $T_x$,
$x\ge 0$. To formalize the argument, note
\[
P(T_{x+y}-T_x>t)=P(\Gamma_x<y,T_{x+y}-T_x>t)
 = P(\Gamma_x<y, X_{T_x+t}-X_{T_x}<y-\Gamma_x).
\]
Since $\Gamma_x=X_{T_x}-x$ is independent of $X_{T_x+t}-X_{T_x}$ and
$X_{T_x+t}-X_{T_x}\stackrel{d}{=}\widetilde X_t$ it follows that
\[
P(T_{x+y}-T_x>t)=P(X_0<y,X_{T_x+t}-X_{T_x}<y-X_0)=P(X_t<y)=P(T_y>t).
\]
\hfill $\Box$

\begin{lemma}\label{lemma_mgf}  For $u>0$ and $\theta<\Phi(u)$, 
\begin{equation} \label{mgf_T}
\int_0^\infty ue^{-ux} E(e^{\theta T_x})\,dx
  = 1 + {\theta\over \Phi(u)-\theta}\, {\Phi(u)\over\mu u}.
\end{equation}
Also, for $u>0$ and $\theta>-\mu u$,
\begin{eqnarray}\label{mgf_Tcenter}
\int_0^\infty ue^{-ux}\,E(e^{\theta (T_x-x/\mu)}-1)\,dx
={\theta^2 \over (\mu u+\theta)^2}\left[ 
 {\mu u \over \Phi(u+\theta/\mu)-\theta}-1\right]
\end{eqnarray}
and
\begin{equation}\label{mgf_Ttildecenter}
\int_0^\infty ue^{-ux}\,E(e^{\theta (\widetilde T_x-x/\mu)}-1)\,dx
={\theta \over \mu u+\theta}\left[ 
 {\mu u \over \Phi(u+\theta/\mu)-\theta}-1\right].
\end{equation}
\end{lemma}
\textbf{Proof}.  Relation (\ref{mgf_T}) follows by letting $v\to 0$ in 
Lemma \ref{jointgammaT} and using that $\Phi(v)/v\to\mu$ in this limit.

The remaining calculations, involving the random variables $T_x-x/\mu$
and $\widetilde T_x-x/\mu$, follow from (\ref{mgf_T}) and the
analogous expression
\[
\int_0^\infty ue^{-ux}\,E(e^{\theta \widetilde T_x})\,dx
={\Phi(u)\over \Phi(u)-\theta},
\]
where we note $\phi(u+\theta/\mu)<\mu u+\theta$ for all $\theta$ such
that $u+\theta/\mu>0$. \hfill $\Box$

\begin{lemma}\label{lemmaincreasing}
The function
\begin{eqnarray*}
E(e^{\theta(T_x-x/\mu)}-1),\quad x\ge 0,
\end{eqnarray*}
is nonnegative for any real parameter $\theta$ and differentiable and
nondecreasing with respect to the variable $x$. The derivative with
respect to $x$ is given by
\[
{d\over dx}E(e^{\theta(T_x-x/\mu)}-1)
=\theta e^{-\theta x/\mu}E(e^{\theta \widetilde T_x}-e^{\theta T_x})/\mu
\ge 0.
\]
\end{lemma}
\textbf{Proof.}  The nonnegativity follows from Jensen's
inequality. It follows from (\ref{mgf_Tcenter}),
(\ref{mgf_Ttildecenter}) and the uniqueness property of Laplace
transforms that $E(e^{\theta(T_x-x/\mu)}-1)$ is obtained as the
convolution of $E(e^{\theta (\widetilde T_x-x/\mu)}-1)$ with the
exponential $e^{-\theta x/\mu}$. Hence
\[ 
E(e^{\theta(T_x-x/\mu)}-1)
={\theta\over\mu}
\int_0^x e^{-\theta(x-y)/\mu}E(e^{\theta (\widetilde T_y-y/\mu)}-1)\,dy.
\]
The left hand side is differentiable in $x$ with derivative
\begin{eqnarray*}
&&{d\over dx}E(e^{\theta(T_x-x/\mu)}-1)
   =-{\theta\over\mu}E(e^{\theta(T_x-x/\mu)}-1)
       +{\theta\over\mu} E(e^{\theta (\widetilde T_x-x/\mu)}-1)\\
&&\qquad  ={\theta\over\mu}e^{-\theta x/\mu}
         E(e^{\theta \widetilde T_x}-e^{\theta T_x}).
\end{eqnarray*}
Now we observe that the processes $T_x$ and $\widetilde T_x$ can be
constructed on the same probability space by a shift of size
$X_0$ so that $T$ is a copy of $\widetilde T$ with the first point of
increase in $X_0$ rather than in 0. In particular $P(\widetilde T_x\ge
T_x)=1$. Hence $\theta E(e^{\theta \widetilde T_x}-e^{\theta T_x})\ge 0$
for any $\theta$. 

\hfill $\Box$ 

\begin{lemma}\label{lemmaTtildemean}
For $x>0$,
\begin{itemize}
\item[i)]
${\displaystyle x/\mu\le E(\widetilde T_x)\le {e\over \Phi(1/x)}}$,
\item[ii)]
${\displaystyle E(\widetilde T_x)\le {e^2(e-1)^{-1}\over \nu(x,\infty)}}$,
\item[iii)]
${\displaystyle{d\over dx}
{\rm Var}(T_x)={2\over\mu}E(\widetilde T_x-x/\mu)\ge 0}$, 
\item[iv)]
${\displaystyle {\rm Var}(T_x)\le {2e\over\mu}\int_0^x \Phi(1/y)^{-1}\,dy}-(x/\mu)^2$.
\end{itemize}
\end{lemma}
\textbf{Proof.}  For i), it was noticed in the proof of Lemma
\ref{lemmaincreasing} that the processes $T_x$ and $\widetilde T_x$
could be constructed such that $\widetilde T_x\ge T_x$ almost
surely. Hence $E(\widetilde T_x)\ge E(T_x)=x/\mu$. Moreover,  
\[
E(\widetilde T_x)=\int_0^\infty P(\widetilde X_t\le x)\,dt
 \le \int_0^\infty e E(e^{-\widetilde X_t/x})\,dt 
  =e\int_0^\infty  e^{-t\Phi(1/x)}\,dt=e/\Phi(1/x).
\]
Inequality ii) follows from
\[
\Phi(1/x)={1\over x}\int_0^\infty e^{-u/x}\nu(u,\infty)\,du
\ge {1\over x}\int_0^x e^{-u/x}\nu(u,\infty)\,du
\ge \nu(x,\infty) (1-e^{-1}).
\]
To prove iii) and iv), differentiate twice with respect to $\theta$ in
(\ref{mgf_Tcenter}) to obtain
\begin{equation}\label{Laplacevariance}
\int_0^\infty u e^{-ux} {\rm Var}(T_x)\,dx
 = {2\over(\mu u)^2}\Big({\mu u \over\Phi(u)}-1\Big).
\end{equation}
Similarly, using (\ref{mgf_Ttildecenter}),
\[
\int_0^\infty u e^{-ux} E(\widetilde T_x-x/\mu)\,dx
 = {1\over \mu u}\Big({\mu u \over\Phi(u)}-1\Big),
\]
hence by partial integration
\[
\int_0^\infty u e^{-ux} \int_0^x E(\widetilde T_y-y/\mu)\,dy\,dx
 = {1\over \mu u^2}\Big({\mu u \over\Phi(u)}-1\Big).
\]
By identification of the Laplace transforms,
\[
{\rm Var}(T_x)={2\over\mu}\int_0^x E(\widetilde T_y-y/\mu)\,dy.
\]
The two inequalities in (i) now imply iii) and iv).

\hfill $\Box$

\subsection{The marginal distribution under scaling}

We will need the weak law of large numbers and an elementary renewal
type theorem for $\widetilde T_x$. Such results are well-known. The
first property below follows from the law of large numbers
for $\widetilde X_t$. The second from the formula
\[
E(\widetilde T_{ax})/a=\int_0^\infty P(\widetilde X_{at}/a\le x)\,dt
  \to x/\mu,
\]
\begin{lemma} \label{lemmaTtildeLLN}
As $a\to\infty$, we have
\begin{itemize}
\item[i)]
${\displaystyle{1\over a}\widetilde T_{ax}\to {x\over\mu}\quad\mbox{in distribution}}$
\item[ii)]
${\displaystyle {1\over a}E(\widetilde T_{ax})\to {x\over\mu}}$. 
\end{itemize}
\end{lemma}

We are now prepared to prove a limit property of the centered variable
$T_x-x/\mu$ under scaling, which is crucial for the distributional
convergence in Theorem \ref{thmICR}.
 
\begin{lemma} \label{lemmaTminusTtilde}
If the sequence $a=a_m$ is such that (\ref{ICR}) holds for some $c>0$,
then as $m\to\infty$,
\begin{equation}\label{TminusTtildelimit}
m \,E(e^{\theta(\widetilde T_{ax}-ax/\mu)/a}-e^{\theta(T_{ax}-ax/\mu)/a})
\to {c^\beta\over\mu\beta}\int_0^x \theta e^{-\theta t/\mu}t^{-\beta}\,dt
\end{equation}
and
\begin{equation}\label{derivativelimit}
m\,{d\over dx}E(e^{\theta(T_{ax}-ax/\mu)/a}-1)
\to{c^\beta\over\mu^2\beta}\int_0^x \theta^2 e^{-\theta t/\mu}t^{-\beta}\,dt.
\end{equation}
\end{lemma}
\textbf{Proof.}  It is enough to prove (\ref{TminusTtildelimit}) since
(\ref{derivativelimit}) then follows directly from Lemma
\ref{lemmaincreasing}.

Recall from (\ref{X_init}) the relation $P(X_0\le
x)={1\over\mu}\int_0^x\nu(u,\infty)\,dy$ where we use the notation
$\nu(y,\infty)=\int_y^\infty \nu(dv)$. For fixed $x$ condition on
$X_0$ to get
\[
P(T_x<t<\widetilde T_x)=P(X_0>x) P(t<\widetilde T_x)
   +{1\over\mu}\int_0^x P(\widetilde T_{x-y}<t<\widetilde T_x)
         \nu(y,\infty)\,dy.
\]
Multiply this identity by $\theta e^{\theta
t}$ and integrate over $t\ge 0$ to obtain
\begin{eqnarray*}
E(e^{\theta \widetilde T_x}-e^{\theta T_x}) =P(X_0>x)E(e^{\theta
\widetilde T_x}-1)+{1\over\mu}\int_0^x E(e^{\theta\widetilde T_x}-e^{\theta
\widetilde T_{x-y}})\,\nu(y,\infty)\,dy.
\end{eqnarray*}
Hence 
\begin{eqnarray}\label{TminusTtilde}\nonumber
\lefteqn{m \,E(e^{\theta \widetilde T_{ax}/a} 
                           -e^{\theta T_{ax}/a})
 = mP(X_0>ax) E(e^{\theta\widetilde T_{ax}/a}-1)}\\
&&\qquad +{1\over\mu}
 \int_0^x E(e^{\theta\widetilde T_{ax}/a}-e^{\theta
          \widetilde T_{a(x-y)}/a})\,am\nu(ay,\infty)\,dy.
\end{eqnarray}
By (\ref{assumpbeta}),    
\[
{1\over \mu} am\,\nu(ay,\infty)\to c^\beta y^{-1-\beta}.
\]
By (\ref{assumpbeta}) and (\ref{ICR}), and using the direct half of
Karamata's theorem,
\[
 mP(X_0>ax) \to \beta^{-1}c^\beta x^{-\beta},
\]
cf.\ Bingham {\it et al.}\ (1987) Thm.\ 1.5.11 ii) (using in their
notation $f(x)=\nu(x,\infty)$, $\rho=-(1+\beta)$, $\sigma=0$).  If we
assume for the moment that the order can be interchanged in which we
integrate over $y$ and take the limit $m,a\to\infty$, then applying the
above asymptotic results as well as Lemma \ref{lemmaTtildeLLN} i), 
\begin{eqnarray*}
\lefteqn{m \,E(e^{\theta\widetilde T_{ax}/a}
                           -e^{\theta T_{ax}/a})}\\
&&\to \beta^{-1}c^\beta x^{-\beta} (e^{\theta x/\mu}-1)
+\int_0^x  (e^{\theta x/\mu}-e^{\theta(x-y)/\mu})c^\beta y^{-1-\beta}\,dy\\
&&=e^{\theta x/\mu} \,{c^\beta\over\mu\beta}
  \int_0^x \theta e^{-\theta t/\mu}t^{-\beta}\,dt,
\end{eqnarray*}
which is the desired relation (\ref{TminusTtildelimit}).  
In the remaining part of the proof we verify the validity of 
this limit operation by deriving an upper bound for the integrand
$E(e^{\theta\widetilde T_{ax}/a}-e^{\theta \widetilde
T_{a(x-y)}/a})\,am\nu(ay,\infty)$ in (\ref{TminusTtilde}), which is
$dy$-integrable over $(0,x]$. 

Using 
\[
\big|E(e^{\theta\widetilde T_x}-e^{\theta\widetilde
T_{x-y}})\big|\le |\theta| 
E\big[(e^{\theta\widetilde T_x}\vee 1)\,
|\widetilde T_x-\widetilde T_{x-y}|\big]
\]
and Hölder's inequality we have, for each integer $k\ge 2$,  
\begin{eqnarray}\label{Ttildediff}
\big|E(e^{\theta\widetilde T_x}-e^{\theta \widetilde
T_{x-y}})\big| \le |\theta| 
E\big[(e^{\theta\widetilde T_x}\vee 1)^{k/(k-1)}\big]^{1-1/k}
E[|\widetilde T_x-\widetilde T_{x-y}|^k]^{1/k}. 
\end{eqnarray}
Now, 
\begin{eqnarray*}
E\big[|\widetilde T_x-\widetilde T_{x-y}|^k\big] &=& E\int_0^\infty\dots
 \int_0^\infty 1_{\{\widetilde T_{x-y}<t_1,\dots,t_k<\widetilde T_x\}}
 \,dt_1\dots dt_k\\
&=& k! \int\dots
 \int_{t_1<\dots<t_k} P(\widetilde T_{x-y}<t_1,\dots,t_k<\widetilde T_x)
 \,dt_1\dots dt_k\\
&=& k! \int\dots
 \int_{t_1<\dots<t_k} P(x-y<\widetilde X_{t_1}<\dots<\widetilde
 X_{t_k} < x)\,dt_1\dots dt_k.
\end{eqnarray*}
For the event $x-y<\widetilde X_{t_1}<\dots<\widetilde X_{t_k} <x$ to
occur it is necessary, in addition to $X_{t_1}\le x$, that all increments
$\widetilde X_{t_j}-\widetilde X_{t_{j-1}}$, $2\le j\le k$ are less
than $y$ in size. Hence the right hand side is at most
\begin{eqnarray*}
k! \int\dots \int_{t_1<\dots<t_k} P(\widetilde X_{t_1}<x,
\widetilde X_{t_j}-\widetilde X_{t_{j-1}}< y,\, 2\le j\le k)\,dt_1\dots dt_k,
\end{eqnarray*}
which equals
\begin{eqnarray*}  
&&k! \int_0^\infty dt_1 P(\widetilde X_{t_1}<x) 
\int_{t_1}^\infty dt_2 P(\widetilde X_{t_2-t_1}<y)\dots 
\int_{t_{k-1}}^\infty dt_k P(\widetilde X_{t_k-t_{k-1}}<y)\\
&&\quad = k! E(\widetilde T_x) E(\widetilde T_y)^{k-1},
\end{eqnarray*}
since the increments of $X(t)$ are independent and stationary.
By (\ref{Ttildediff}), 
\begin{eqnarray*}
\big|E(e^{\theta\widetilde T_{ax}/a}-e^{\theta \widetilde
T_{a(x-y)}/a})\big|^k\le k! |\theta|^k
E\big[(e^{\theta\widetilde T_{ax}/a}\vee 1)^{k/(k-1)}\big]^{k-1}
E(\widetilde T_{ax}/a) E(\widetilde T_{ay}/a)^{k-1}.
\end{eqnarray*}
By Lemma \ref{lemmaTtildeLLN} we may assume
\[
E\big[(e^{\theta\widetilde T_{ax}/a}\vee 1)^{k/(k-1)}\big]^{k-1}
E(\widetilde T_{ax}/a) \le 2 (e^{\theta x/\mu}\vee 1)^k (x/\mu), 
\]
and thus
\begin{eqnarray*}
\big|E(e^{\theta\widetilde T_{ax}/a}-e^{\theta \widetilde
T_{a(x-y)}/a})\big|
\le C_{\theta,k}(x)
\, E(\widetilde T_{ay}/a)^{1-1/k}
\end{eqnarray*}
for $a\ge a_0$ and sufficiently large $a_0$, with $C_{\theta,k}(x)=
|\theta|(2\,k!)^{1/k} \,(e^{\theta x/\mu}\vee 1)\,(x/\mu)^{1/k}$.
For the integrand in (\ref{TminusTtilde}) we have obtained
\begin{eqnarray}\label{domconv}
\lefteqn{\big|E(e^{\theta\widetilde T_{ax}/a}-e^{\theta \widetilde
T_{a(x-y)}/a})\big|\, am\nu(ay,\infty)}\\
&&\le C_{\theta,k}(x)\, E(\widetilde T_{ay}/a)^{1-1/k}\,
am\nu(ay,\infty),\quad 0<y\le x, \;a\ge a_0. \nonumber
\end{eqnarray}
We split the further task of estimating the right hand side in the
above expression in the two cases $ay>a_0$ and $ay\le a_0$.

By (\ref{muassumption}) and (\ref{lowerindex}), there exist a constant
$C_1$ such that for any $q<\sigma$ we have $\Phi(\lambda)\ge
C_1(\lambda\wedge \lambda^q)$, $\lambda>0$. The lower bound in
(\ref{lowerindex}) ensures, moreover, that we may take $q$ such that
$\beta<q<\sigma$. In combination with Lemma \ref{lemmaTtildemean} i),
this yields, for such $q$, $E(\widetilde T_x)\le C_2(x\vee x^q)$. Thus,
\[
E(\widetilde T_{ay}/a)\le  C_2(y\vee (y^q/a_0^{1-q}))\le C_3\, y^q,
\quad 0\le y\le x, \quad a\ge a_0.
\]
Furthermore, since the function $\nu(x,\infty)$ is regularly varying
at infinity with index $-(1+\beta)$, we have for $ay>a_0$ and
$\epsilon>0$ the Potter type bound
\[
am\nu(ay,\infty)\le C_4 \,y^{-1-\beta} \max(y^\epsilon,y^{-\epsilon})
\]
(Bingham {\it et al.}\ (1987), Ch.\ 1.5). Thus, for some constant $C$,
\[
E(\widetilde T_{ay}/a)^{1-1/k}\,am\nu(ay,\infty)\le 
   C\,y^{q(1-1/k)}\, y^{-1-\beta-\epsilon}. 
\]
Since $q>\beta$ we may take $k$ so large that $q(1-1/k)>\beta$ and
then $\epsilon$ so small that $\epsilon<q(1-1/k)-\beta$ to obtain a
dominating function for the integrand in (\ref{domconv}) which is
integrable in $y$ over $[0,x]$.

For the remaining case $ay\le a_0$, Lemma \ref{lemmaTtildemean}
i) implies
\[
E(\widetilde T_{ay}/a)^{1-1/k}\,am\nu(ay,\infty)
\le  (e^2/(e-1))^{1-1/k}\, m\,a^{1/k} \nu(ay,\infty)^{1/k}.
\]
Using a property of slowly varying functions (Bingham {\it et al.}\
(1987), Prop 1.3.6), for any $\epsilon>0$, $L(a)a^\epsilon\to
\infty$ as $a\to\infty$. Hence we may assume $a^{-\epsilon}\le L(a)$.
Also,
\[
\nu(ay,\infty)\le {1\over ay}\int_{ay}^\infty u\,\nu(du)\le {\mu \over ay}.
\]
Thus, using (\ref{ICR}), 
\[
m\,a^{1/k} \nu(ay,\infty)^{1/k}\le {m L(a)\over a^\beta}
a^{\epsilon+\beta+1/k}\,(\mu/ay)^{1/k}
\le 2c^\beta \mu^{1+1/k} a^{\epsilon+\beta} y^{-1/k}.
\]
Now apply $a\le a_0/y$ to obtain from (\ref{domconv}) a constant $C$
for which 
\begin{eqnarray*}
\big|E(e^{\theta\widetilde T_{ax}/a}-e^{\theta \widetilde
T_{a(x-y)}/a})\big|\, am\nu(ay,\infty)
\le C\, a_0^{\epsilon+\beta}\, y^{-\beta-1/k-\epsilon}.
\end{eqnarray*}
This is again integrable if we make the same choise of $k$ and
$\epsilon$ as above. This concludes the proof that the limit in
(\ref{TminusTtilde}) can be carried out under the integral sign and
hence the proof of the lemma. 

\hfill $\Box$

\begin{lemma}\label{lemmamargscaling}
For $a=a_m$ such that (\ref{ICR}) holds for some $c>0$, 
\[ 
mE(e^{\theta(T_{ax}-ax/\mu)/a}-1)\to 
  {c^\beta\over\beta\mu^2}\int_0^x \int_0^y 
   \theta^2 e^{-\theta s/\mu}s^{-\beta}\,dsdy,\quad m\to\infty.
\]
\end{lemma}
\textbf{Proof}.  By Lemma \ref{lemmaincreasing},
$mE(e^{\theta(T_{ax}-ax/\mu)/a}-1)$ is nonnegative and increasing in
$x$. The limit function on the right hand side is also nonnegative and
increasing. Hence the lemma follows from weak convergence of measures
if we can prove
\begin{equation}\label{Laplaceconv} 
\int_0^\infty  e^{-ux} {d\over dx}mE(e^{\theta(T_{ax}-ax/\mu)/a}-1)\,dx
\to 
\int_0^\infty e^{-ux}\left({c^\beta\over\mu^2\beta}\int_0^x 
    \theta^2 e^{-\theta s/\mu}s^{-\beta}\,ds\right)dx.
\end{equation}
To find the Laplace transform on the right hand side note that
\[
{\theta^2\Gamma(1-\beta)\over \beta (u-\theta)^{1-\beta}}
= {\theta^2\over\beta} \int_0^\infty e^{-ux} e^{\theta x}
x^{-\beta}\,dx,\quad  \theta <u.
\]
Multiplication of the transform by $1/u$ corresponds to
integration of $e^{\theta x} x^{-\beta}$. Hence 
\[
\int_0^\infty e^{-ux}\left({1\over\beta}\int_0^x 
        \theta^2 e^{\theta s}s^{-\beta}\,ds\right)dx
= {\Gamma(1-\beta)\theta^2\over\beta u (u-\theta)^{1-\beta}},\quad
        \theta < u,
\]
and hence (\ref{Laplaceconv}) is equivalent to  
\begin{equation}\label{ICR_limit}
\int_0^\infty ue^{-ux} mE(e^{\theta(T_{ax}-ax/\mu)/a}-1)\,dx
\to {\Gamma(1-\beta)c^\beta\theta^2\over\beta u (u+\theta/\mu)^{1-\beta}\mu^2},\quad
        \theta > -\mu u.
\end{equation}
To help analyze the Laplace transform in (\ref{ICR_limit})
we introduce the additional notation 
\[
I(u)=\mu u-\Phi(u)=\int_0^\infty (e^{-ux}-1+ux)\,\nu(dx)
\ge 0.
\]
Writing $I(u)=u^2\int_0^\infty e^{-ux}U(x)\,dx$ with
$U(x)=\int_x^\infty \nu(y,\infty)\,dy$, it follows from Karamata's
Tauberian Theorem (Thm.\ 1.7.6 in Bingham {\it et al.}\ (1987)) that 
\begin{equation}\label{I_karamata}
aI(u/a)\sim {\Gamma(1-\beta)L(a/u)u^{1+\beta}\over \beta a^\beta}, 
\quad  a\to\infty 
\end{equation}
Relation (\ref{mgf_Tcenter}) of Lemma \ref{lemma_mgf} now shows
\begin{eqnarray*}
\lefteqn{\int_0^\infty ue^{-ux} mE(e^{\theta (T_{ax}-ax/\mu)/a}-1)\,dx}\\
&&={m\theta^2 \over (\mu u+\theta)^2}\, {aI((u+\theta/\mu)/a)\over \mu
u-aI((u+\theta/\mu)/a)} \\ 
&&\sim {m L(a)\over a^\beta}
\left(u+\theta/\mu\right)^{-(1-\beta)}\, {\Gamma(1-\beta)\theta^2\over
\beta\mu^3 u} 
\\
&&\sim {\Gamma(1-\beta)c^\beta\theta^2\over\beta u(u+\theta/\mu)^{1-\beta}
\mu^2},  \quad \theta>-\mu u,
\end{eqnarray*}
which proves (\ref{ICR_limit}) and hence the lemma. 
\hfill $\Box$

We are now able to conclude convergence of the marginal
distributions. 
\begin{lemma} \label{lemmamargdistr}
Under the assumptions of Theorem 2, for any $x\ge 0$
\[
{1\over a_m}\sum_{i=1}^m (T^i_{a_m x}-{1\over\mu} a_m x)
\stackrel{d}{\to} -{1\over\mu} c\,Y_\beta(x/c)
\]
\end{lemma}
\textbf{Proof}.  
Writing
\[
\Lambda^{(m)}(\theta;x)=mE(e^{\theta(T_{ax}-ax/\mu)/a}-1),
\]
Lemma \ref{lemmamargscaling} shows that 
\begin{eqnarray*}
\lefteqn{\log E\exp\Big\{\theta  
 {1\over a_m}\sum_{k=1}^m (T^{(k)}_{a_m x}-{1\over\mu}a_m x)
       \Big\}\,dx }\\
&&=\log\Big(1+{1\over m}\Lambda^{(m)}(\theta;x)\Big)^m
\to {c^\beta\over\beta\mu^2}\int_0^x \int_0^y 
        \theta^2 e^{-\theta s/\mu}s^{-\beta}\,dsdy.
\end{eqnarray*}
This proves the lemma since the limit process $Y_\beta$ has the property
\[
\log E(e^{\theta Y_\beta(x)})={1\over\beta}\int_0^x \int_0^y 
        \theta^2 e^{\theta s}s^{-\beta}\,dsdy
\]
and so, as noticed in Lemma \ref{lemmaconsistency},
\[
\log E(e^{-\theta cY_\beta(x/c)/\mu})
 ={c^\beta\over\beta\mu^2}\int_0^x \int_0^y 
        \theta^2 e^{-\theta s/\mu}s^{-\beta}\,dsdy.
\] 
\hfill$\Box$

\section{Multivariate distributions}

The proofs of convergence of the finite-dimensional distributions are
based on the following recursive equations for moment generating
functions.

\begin{lemma}  \label{mgfprop}
Fix $n\ge 2$ and a sequence of time points $0\le x_1\le \dots\le
x_n$. The moment generating function of the finite-dimensional
distributions of the stationary inverse Lévy subordinator process $\{T_x\}$
satisfies the recurrence relation
\begin{eqnarray}
\lefteqn{E\exp\Big\{\sum_{i=1}^n \theta_i T_{x_i}\Big\} 
 = E\exp\Big\{\sum_{i=2}^n \theta_i T_{x_i}\Big\}} \nonumber \\ 
&& + {\theta_1\over\sum_{i=1}^n \theta_i}{\displaystyle \int_0^{x_1}} 
E\Big[\exp\Big\{\sum_{i=2}^n \theta_i \widetilde T_{x_i-x}\Big\}\Big]
\,d_xE\Big[\exp\Big\{T_x\,\sum_{i=1}^n \theta_i\Big\}\Big], 
\label{mgfstat}
\end{eqnarray}
where $\widetilde T_x$ is the corresponding pure inverse Lévy process.
Moreover,
\begin{eqnarray}
\lefteqn{E\exp\Big\{\sum_{i=1}^n \theta_i \widetilde T_{x_i}\Big\} 
 = E\exp\Big\{\sum_{i=2}^n \theta_i \widetilde T_{x_i}\Big\}} \nonumber \\ 
&& + {\theta_1\over\sum_{i=1}^n \theta_i}{\displaystyle \int_0^{x_1}} 
E\Big[\exp\Big\{\sum_{i=2}^n \theta_i \widetilde T_{x_i-x}\Big\}\Big]
\,d_xE\Big[\exp\Big\{\widetilde T_x\,\sum_{i=1}^n \theta_i\Big\}\Big], 
\label{mgfpure}
\end{eqnarray}
\end{lemma}

\noindent \textbf{Proof}. 
We have
\[
E\exp\Big\{\sum_{i=1}^n \theta_i T_{x_i}\Big\}-
               E\exp\Big\{\sum_{i=2}^n \theta_i T_{x_i}\Big\} 
=E\Big[\exp\Big\{\sum_{i=2}^n \theta_i T_{x_i}\Big\}
    \Big(e^{\theta_1 T_{x_1}}-1\Big)\Big].
\]
Since 
\[
e^{\theta_1 T_{x_1}}-1
 =\int_0^\infty 1_{\{u\le T_{x_1}\}}\theta_1 e^{\theta_1 u}\,du
 =\int_0^\infty 1_{\{X_u\le x_1\}}\theta_1 e^{\theta_1 u}\,du,
\]
it follows that 
\begin{eqnarray*}
\lefteqn{E\Big[ \exp\Big\{\sum_{i=2}^n \theta_i T_{x_i}\Big\}
           \Big(e^{\theta_1 T_{x_1}}-1\Big)\Big]}\\
&&=E\Big[\int_0^\infty 1_{\{X_u\le x_1\}} 
    \exp\Big\{\sum_{i=2}^n \theta_i T_{x_i}\Big\} 
                                \theta_1 e^{\theta_1 u}\,du \Big]\\
  && = E\Big[\int_0^\infty 1_{\{X_u\le x_1\}} 
         \exp\Big\{\sum_{i=2}^n \theta_i(T_{x_i}-T_{X_u})\Big\}\, 
          \theta_1 \exp\Big\{(\sum_{i=1}^n \theta_i) T_{X_u}\Big\}\,du \Big].
\end{eqnarray*}
Here, $T_{X_u}=u$. For any $u>0$ and $i\ge 2$, on the set
$\{X_u\le x_1\}$ we have
\[
\{T_{x_i}-T_{X_u}\le t\}=\{T_{x_i}\le u+t\}=\{X_{u+t}> x_i\}.
\]
Since $\{X_t\}$ has independent increments the rightmost event has the
same probability as 
\[ 
\{X_u+\widetilde X_t> x_i\}=\{\widetilde T_{x_i-X_u}\le t\},
\]
where $X_u\le x_1$ is assumed independent of $\widetilde X_t$. Thus,
on $\{X_u\le x_1\}$ the increment $T_{x_i}-T_{X_u}$ has the same
distribution as $\widetilde T_{x_i-X_u}$. It follows that
\begin{eqnarray*}
\lefteqn{E\exp\Big\{\sum_{i=1}^n \theta_i T_{x_i}\Big\} 
 - E\exp\Big\{\sum_{i=2}^n \theta_i T_{x_i}\Big\}} \nonumber \\ 
  &&=\theta_1  E\Big[\int_0^\infty 1_{\{X_u\le x_1\}} 
   E\Big[\exp\Big\{\sum_{i=2}^n \theta_i\widetilde T_{x_i-X_u}\Big\}
                                                |X_u\Big] 
      \exp\Big\{(\sum_{i=1}^n \theta_i)u\Big\} \,  du\Big]\\
  &&=\theta_1  E\Big[\int_0^{x_1} 
     E\Big[\exp\Big\{\sum_{i=2}^n \theta_i \widetilde T_{x_i-x}\Big\}\Big] 
      \exp\Big\{(\sum_{i=1}^n \theta_i)T_{x}\Big\}\, dT_x\Big],
\end{eqnarray*}
where the integration after variable substitution $x=X_u$ is with
respect to the increasing function of bounded variation $\{T_x, x\ge
0\}$. (Intuitively, the time-change $X_u$ picks out the rightmost
point of each flat piece of $T_x$.)  Moreover, if we change to the
measure
\begin{eqnarray*}
d_x\Big(\exp\Big\{T_x\,\sum_{i=1}^n \theta_i\Big\}\Big)
  =\Big(\sum_{i=1}^n \theta_i\Big) 
          \exp\Big\{T_x\,\sum_{i=1}^n \theta_i \Big\}\,dT_x
\end{eqnarray*}
we obtain
\begin{eqnarray*}
\lefteqn{E\exp\Big\{\sum_{i=1}^n \theta_i T_{x_i}\Big\} 
 - E\exp\Big\{\sum_{i=2}^n \theta_i T_{x_i}\Big\}} \nonumber \\ 
  &&\hspace{-2mm}={\theta_1 \over \sum_{i=1}^n \theta_i}
  \, E\Big[\int_0^{x_1} E\Big[\exp\Big\{\sum_{i=2}^n \theta_i 
                                         \widetilde T_{x_i-x}\Big\}\Big]\, 
         d_x\Big(\exp\Big\{T_x\,\sum_{i=1}^n\theta_i\Big\}\Big)\Big]\\
  &&\hspace{-2mm}={\theta_1 \over \sum_{i=1}^n \theta_i}\, 
  \int_0^{x_1}E\Big[\exp\Big\{\sum_{i=2}^n \theta_i 
                                    \widetilde T_{x_i-x}\Big\}\Big] \, 
         d_x E\Big[\exp\Big\{T_x\,\sum_{i=1}^n\theta_i\Big\}\Big],
\vspace{-25pt}
\end{eqnarray*}
which is (\ref{mgfstat}). Start with $\widetilde T$ rather than $T$ to get
(\ref{mgfpure}).\hfill $\Box$

For $n\ge 1$ and $1\le k\le n$, put
$\bar\theta_{k,n}=(\theta_k,\dots,\theta_n)$ and $\bar
x_{k,n}=(x_k,\dots,x_n)$, where $0=x_0\le x_1\le \dots\le x_n$ and let
\begin{equation}\label{defPhi}
\Phi_{n-k+1}(\bar\theta_{k,n};\,\bar x_{k,n})
 = E\exp\Big\{\sum_{i=k}^n \theta_i(T_{x_i}-x_i/\mu)\Big\}  
\end{equation}
denote the multivariate moment generating functions for the centered
process $\{T_x-x/\mu\}_{x\ge 0}$.  Here, the subindex $n-k+1$ is the
number of elements of the argument vectors $\bar\theta_{k,n}$, $\bar
x_{k,n}$.  Similarly, let
$\widetilde\Phi_{n-k+1}(\bar\theta_{k,n};\,\bar x_{k,n})$, $1\le k\le
n$, denote the corresponding functions for the pure process
$\{\widetilde T_x-x/\mu\}_{x\ge 0}$.  The subtraction $\bar
x_{k,n}-u=(x_k-u,\dots,x_n-u)$ is interpreted component-wise in the
next statement and in the sequel. 
\begin{lemma} \label{lemmaPhi} The moment generating functions defined in 
(\ref{defPhi}) satisfy the integral equation
\begin{eqnarray*} 
\lefteqn{\Phi_n(\bar\theta_{1,n};\,\bar x_{1,n})
=\Phi_{n-1}(\bar\theta_{2,n};\,\bar x_{2,n})e^{-\theta_1 x_1/\mu}
}\\[2mm] &&+{\theta_1\over \sum_{i=1}^n\theta_i} \int_0^{x_1}
e^{-\theta_1(x_1-x)/\mu} \,\widetilde\Phi_{n-1}(\bar\theta_{2,n};\bar
x_{2,n}-x)\, \Phi_1\Big(\sum_{i=1}^n\theta_i;\,dx\Big)\\
&&+{\theta_1\over \mu} \int_0^{x_1} e^{-\theta_1(x_1-x)/\mu}
\,\widetilde\Phi_{n-1}(\bar\theta_{2,n};\bar x_{2,n}-x)\,
\Phi_1\Big(\sum_{i=1}^n\theta_i;\,x\Big)\,dx.
\end{eqnarray*}
\end{lemma}
\textbf{Proof}.  By Lemma \ref{mgfprop},
\begin{eqnarray*}\label{Phieqn}
\lefteqn{\Phi_n(\bar\theta_{1,n};\,\bar x_{1,n}) 
=\Phi_{n-1}(\bar\theta_{2,n};\,\bar x_{2,n})e^{-\theta_1 x_1/\mu} 
 +{\theta_1\over \sum_{i=1}^n\theta_i} }\\[2mm] 
&&\times\int_0^{x_1} e^{-\theta_1(x_1-x)/\mu} \nonumber
   \widetilde\Phi_{n-1}(\bar\theta_{2,n};\bar x_{2,n}-x)
    \exp\Big\{\!-\!{x\over\mu}\sum_{i=1}^n\theta_i\Big\}
            \,d_xE\Big[\exp\Big\{T_x\sum_{i=1}^n\theta_i\Big\}\Big], 
\end{eqnarray*}
which, by observing
\begin{eqnarray*}
\lefteqn{ \exp\Big\{\!-\!{x\over\mu}\sum_{i=1}^n\theta_i\Big\}
             \, d_xE\Big[\exp\Big\{T_x\sum_{i=1}^n\theta_i\Big\}\Big]}\\
&&= d_xE\Big[\exp\Big\{(T_x-x/\mu)\sum_{i=1}^n\theta_i\Big\}\Big]
    + {1\over\mu}\sum_{i=1}^n\theta_i\,
        E\Big[\exp\Big\{(T_x-x/\mu)\sum_{i=1}^n\theta_i\Big\}\Big]\,dx\\
&&= \Phi_1\Big(\sum_{i=1}^n\theta_i;\,dx\Big)
    +{1\over\mu}\sum_{i=1}^n\theta_i\,
                 \Phi_1\Big(\sum_{i=1}^n\theta_i;\,x\Big)\,dx,
\end{eqnarray*}
may be rewritten in the form stated in the lemma.
\hfill $\Box$

According to (\ref{defmgf}) we must find the
limits of the scaled function
\[
m(\Phi_n(\bar\theta_{1,n}/b;\,a \bar x_{1,n})-1) 
  = mE\Big[\exp\Big\{\sum_{i=1}^n \theta_i(T_{ax_i}-ax_i/\mu)/b\Big\}-1\Big]  
\]
as $m$, $a$ and $b$ tend to infinity, when $a$ and $b$ satisfy either
(\ref{FCR}) together with (\ref{FCR_b}) or condition (\ref{ICR}).
The first case is {\it FBM scaling} leading to fractional Brownian motion
in the limit, as in Theorem 1, and the second case (with $a=b$) is the {\it
intermediate scaling} studied in Theorem 2.

For $n\ge 1$, $m\ge 1$ and $a,b>0$ we introduce
\begin{eqnarray*}
\Lambda_n^{(m)}(\bar\theta_{1,n};\,\bar x_{1,n}) &=&
 m(\Phi_n(\bar\theta_{1,n}/b;\,a \bar x_{1,n})-1),
\end{eqnarray*}
as well as
\begin{eqnarray*}
\widetilde\Lambda_n^{(m)}(\bar\theta_{1,n};\,\bar x_{1,n}) &=&
 {am\over b}(\widetilde\Phi_n(\bar\theta_{1,n}/b;\,a \bar x_{1,n})-1)
\end{eqnarray*}
and
\begin{eqnarray}\label{defxi}
\Xi^{(m)}_n(\bar\theta_{1,n};\,\bar x_{1,n}) 
&=&\widetilde\Lambda_n^{(m)}(\bar\theta_{1,n};\,\bar x_{1,n}) 
  - \Lambda_n^{(m)}(\bar\theta_{1,n};\,\bar x_{1,n}).
\end{eqnarray}
Our strategy for finding the corresponding limit functions is to
derive for fixed $m$ sequences of integral equations, which are
recursive in $n$. As already pointed out we give the detailed proof
only for Theorem 2. To simplify notation during this analysis we will also
use a special notation for the sums
\[
\eta_n=\sum_{i=1}^n\theta_i,\quad n\ge 1,
\]
which appear frequently as evident in Lemma \ref{lemmaPhi}. 

\subsection{Multivariate distribution under the intermediate scaling}

We study the asymptotic limits of $\Lambda_n^{(m)}$
and $\widetilde\Lambda_n^{(m)}$ as $m\to\infty$ under assumption
(\ref{ICR}). For simplicity the constant in (\ref{ICR}) is set to
$c=1$. The general case $c\not=1$ then follows from Lemma
\ref{lemmaconsistency}. We begin with a system of equations for the
functions $\Xi^{(m)}_n$ defined in (\ref{defxi}), which will be used
to determine corresponding limit functions as $m\to\infty$.

\begin{lemma} \label{lemmaXi_m}
We have
\begin{eqnarray*}\nonumber
\lefteqn{\Xi_n^{(m)}(\bar\theta_{1,n};\,\bar x_{1,n}) 
\;=\; \Xi^{(m)}_1 (\eta_n;\,x_1)
  +e^{-\theta_1x_1/\mu}\,
\Xi^{(m)}_{n-1}(\bar\theta_{2,n};\bar \,x_{2,n})}\\
&& + (\theta_1/\eta_n-1)
    \int_0^{x_1} e^{-\theta_1(x_1-x)/\mu}\,
            \Xi^{(m)}_1 (\eta_n;\,dx)
          +{1\over m}R^{(m)},
\end{eqnarray*}
where
\[
R^{(m)}=\frac{\theta_1}{\eta_n} 
 \int_0^{x_1} e^{-\theta_1(x_1-x)/\mu} 
   \,\widetilde\Lambda^{(m)}_{n-1}(\bar\theta_{2,n};\bar x_{2,n}-x)\,
           \widetilde\Lambda^{(m)}_1(\eta_n;\,dx) 
\]
\end{lemma}
\noindent
\textbf{Proof}. After inserting the scaling parameters $m$ and $a$ into
the equation obtained in Lemma \ref{lemmaPhi} and sorting the terms
appropriately, it is seen that the scaled functions $\Lambda_n^{(m)}$
and $\widetilde\Lambda_n^{(m)}$ satisfy 
\begin{eqnarray*}
\Lambda_n^{(m)}(\bar\theta_{1,n};\,\bar x_{1,n}) 
=\Lambda_{n-1}^{(m)}(\bar\theta_{2,n};\,\bar x_{2,n})e^{-\theta_1x_1/\mu}
 +I_1^{(m)}+I_2^{(m)}+I_3^{(m)}+{1\over m}R_1^{(m)},
\end{eqnarray*}
where
\begin{eqnarray*}
I_1^{(m)}&=&\frac{\theta_1}{\eta_n} 
 \int_0^{x_1} e^{-\theta_1(x_1-x)/\mu} 
           \Lambda^{(m)}_1(\eta_n;\,dx) 
\\
I_2^{(m)}&=& {\theta_1\over \mu}
  \int_0^{x_1} e^{-\theta_1(x_1-x)/\mu}
   \,\widetilde\Lambda^{(m)}_{n-1}(\bar\theta_{2,n};\bar x_{2,n}-x)\,dx\\
I_3^{(m)}&=& {\theta_1\over \mu}
  \int_0^{x_1} e^{-\theta_1(x_1-x)/\mu}
          \Lambda^{(m)}_1(\eta_n;\,x) 
\,dx\\
R_1^{(m)}&=&\frac{\theta_1}{\eta_n} 
 \int_0^{x_1} e^{-\theta_1(x_1-x)/\mu} 
   \,\widetilde\Lambda^{(m)}_{n-1}(\bar\theta_{2,n};\bar x_{2,n}-x)\,
           \Lambda^{(m)}_1(\eta_n;\,dx)
\\
&&+{\theta_1\over \mu}
  \int_0^{x_1} e^{-\theta_1(x_1-x)/\mu}
   \,\widetilde\Lambda^{(m)}_{n-1}(\bar\theta_{2,n};\bar x_{2,n}-x)\,
         \Lambda^{(m)}_1(\eta_n;x)
\,dx
\end{eqnarray*}
By Lemma \ref{lemmaincreasing},
\begin{eqnarray}\nonumber
&&{d\over dx}\Lambda^{(m)}_1(\theta;\,x)=
{d\over dx} mE(e^{\theta(T_{ax}-ax/\mu)/a}-1)\\
&&\qquad  = {\theta\over\mu} m E(e^{\theta(\widetilde T_{ax}-ax/\mu)/a}
-e^{\theta(T_{ax}-x/\mu)/a}) = {\theta\over\mu}\Xi_1^{(m)}(\theta;\, x).
\label{derivativeonedim}
\end{eqnarray}
Thus, 
\[
R_1^{(m)}={\theta_1\over \mu}
  \int_0^{x_1} e^{-\theta_1(x_1-x)/\mu}
   \,\widetilde\Lambda^{(m)}_{n-1}(\bar\theta_{2,n};\bar x_{2,n}-x)\,
         \widetilde\Lambda^{(m)}_1(\eta_n;\,x)
\,dx.
\]
By partial integration,
\[
I_3^{(m)} = \Lambda^{(m)}_1(\eta_n;\,x_1)  
- \int_0^{x_1} e^{-\theta_1(x_1-x)/\mu}
          \Lambda^{(m)}_1(\eta_n;\,dx). 
\]
Hence
\begin{eqnarray} \nonumber
\Lambda_n^{(m)}(\bar\theta_{1,n};\,\bar x_{1,n}) 
&=&  \Lambda_{n-1}^{(m)}(\bar\theta_{2,n};\,\bar x_{2,n})e^{-\theta_1x_1/\mu}
   +\Lambda^{(m)}_1(\eta_n;\,x_1)  
\\ \label{Lambdaeqn}
&& +I_2^{(m)}+I_4^{(m)}+{1\over m}R_1^{(m)},
\end{eqnarray}
where now
\begin{eqnarray*}
I_4^{(m)}&=&  \Big(\frac{\theta_1}{\eta_n}-1\Big)
    \int_0^{x_1} e^{-\theta_1(x_1-x)/\mu} 
       \Lambda^{(m)}_1(\eta_n;\,dx). 
\end{eqnarray*}
Similarly,
\begin{eqnarray}\nonumber
\widetilde\Lambda_n^{(m)}(\bar\theta_{1,n};\,\bar x_{1,n}) 
&=&\widetilde\Lambda^{(m)}_{n-1}(\bar\theta_{2,n};\bar x_{2,n})
   e^{-\theta_1x_1/\mu}
   +\widetilde\Lambda^{(m)}_1(\eta_n;\,x_1) 
\\
&& +I_2^{(m)}+\widetilde I_4^{(m)}+{1\over m}R_2^{(m)}, 
\label{tildeLambdaeqn}
\end{eqnarray}  
with 
\begin{eqnarray*}
\widetilde I_4^{(m)}&=&\Big(\frac{\theta_1}{\eta_n}-1\Big)
    \int_0^{x_1} e^{-\theta_1(x_1-x)/\mu} 
            \widetilde\Lambda^{(m)}_1(\eta_n;\,dx)
\end{eqnarray*}
and 
\begin{eqnarray*}
R_2^{(m)}&=&\frac{\theta_1}{\eta_n}
 \int_0^{x_1} e^{-\theta_1(x_1-x)/\mu} 
   \,\widetilde\Lambda^{(m)}_{n-1}(\bar\theta_{2,n};\bar x_{2,n}-x)\,
           \widetilde\Lambda^{(m)}_1(\eta_n;\,dx) 
\\
&&+{\theta_1\over \mu}
  \int_0^{x_1} e^{-\theta_1(x_1-x)/\mu}
   \,\widetilde\Lambda^{(m)}_{n-1}(\bar\theta_{2,n};\bar x_{2,n}-x)\,
         \widetilde\Lambda^{(m)}_1(\eta_n;\,x)
\,dx.
\end{eqnarray*}
By subtracting (\ref{Lambdaeqn}) from (\ref{tildeLambdaeqn}) and using
$R^{(m)}=R_2^{(m)}-R_1^{(m)}$ we obtain the desired equation for
$\Xi^{(m)}_n$.

\hfill $\Box$

The next result generalizes Lemma \ref{lemmaincreasing} to the
multivariate distributions. 
  
\begin{lemma}\label{lemmaderivative_m} 
For any $0\le s\le x_1$ and $n\ge 1$,
\[ 
-{d\over ds} \Lambda^{(m)}_n(\bar\theta_{1,n};\,\bar x_{1,n}-s)  
  =  \frac{\eta_n}{\mu}\, 
           \Xi^{(m)}_n(\bar\theta_{1,n};\,\bar x_{1,n}-s). 
\]
\end{lemma}   
\textbf{Proof}. For $n=1$ this is (\ref{derivativeonedim}).  For $n\ge
2$ and $0<s<x_1$, using (\ref{Lambdaeqn}),
\begin{eqnarray*}
\lefteqn{{d\over ds}\Lambda^{(m)}_n(\bar\theta_{1,n};\,\bar x_{1,n}-s) 
= {1\over m}{d\over ds} R_1^{(m)}(\bar x_{1,n}-s)
+\frac{\theta_1}{\eta_n} 
{d\over ds}\Lambda^{(m)}_1(\eta_n;\,x_1-s) 
}\\ 
&&+  \Big(\frac{\theta_1}{\eta_n}-1\Big)
   {\theta_1\over\mu} \int_0^{x_1-s} e^{-\theta_1(x_1-s-x)/\mu} 
            \Lambda^{(m)}_1(\eta_n;\,dx) 
\\
&& + e^{-\theta_1(x_1-s)/\mu}
   {d\over ds} \Lambda^{(m)}_{n-1}(\bar\theta_{2,n};\bar x_{2,n}-s)
-{\theta_1\over\mu} e^{-\theta_1(x_1-s)/\mu}\,
   \Xi^{(m)}_{n-1}(\bar\theta_{2,n};\bar x_{2,n}-s)
\end{eqnarray*}
and hence by (\ref{derivativeonedim}),
\begin{eqnarray*}
\lefteqn{-{d\over ds}\Lambda^{(m)}_n(\bar\theta_{1,n};\,\bar x_{1,n}-s) 
= -{1\over m}{d\over ds} R_1^{(m)}(\bar x_{1,n}-s)
+{\theta_1\over\mu}\Xi^{(m)}_1(\eta_n;\,x_1-s) 
}\\ 
&&+ e^{-\theta_1(x_1-s)/\mu} \frac{1}{\mu}(\eta_n-\theta_1) 
  {\theta_1\over\mu}\int_0^{x_1-s} e^{\theta_1x/\mu}\, 
            \Xi^{(m)}_1(\eta_n;\,x) 
\,dx\\
&&- e^{-\theta_1(x_1-s)/\mu}\left(
{d\over ds} \Lambda^{(m)}_{n-1}(\bar\theta_{2,n};\bar x_{2,n}-s)
- {\theta_1\over\mu}\Xi_{n-1}(\bar\theta_{2,n};\bar x_{2,n}-s)\right).
\end{eqnarray*}
Here, 
\begin{eqnarray*}
\lefteqn{{\theta_1\over\mu}\int_0^{x_1-s} e^{\theta_1x/\mu}\, 
            \Xi^{(m)}_1(\eta_n;\,x) 
\,dx}\\
&&\!\! = e^{\theta_1(x_1-s)/\mu} \Xi^{(m)}_1(\eta_n;\,x_1-s) 
     -\int_0^{x_1-s} e^{\theta_1 x/\mu}\,\Xi^{(m)}_1(\eta_n;\,dx)
\end{eqnarray*}
so
\begin{eqnarray*}  
&&\hspace{-3mm} -{d\over ds}\Lambda^{(m)}_n(\bar\theta_{1,n};\,\bar x_{1,n}-s) 
= -{1\over m}{d\over ds} R_1^{(m)}(\bar x_{1,n}-s)
+\frac{\eta_n}{\mu}\, 
 \Xi^{(m)}_1(\eta_n;\,x_1-s) 
\\ 
&& \qquad\qquad -\frac{1}{\mu}(\eta_n-\theta_1) 
     \int_0^{x_1-s} e^{-\theta_1(x_1-x-s)/\mu}\, 
            \Xi^{(m)}_1(\eta_n;\,dx) 
\\
&&\qquad\qquad -e^{-\theta_1(x_1-s)/\mu}\left(
   {d\over ds} \Lambda^{(m)}_{n-1}(\bar\theta_{2,n};\bar x_{2,n}-s)
 -{\theta_1\over\mu}\Xi^{(m)}_{n-1}(\bar\theta_{2,n};\bar x_{2,n}-s)\right).
\end{eqnarray*}
By replacing the integral term using Lemma \ref{lemmaXi_m} 
this implies
\begin{eqnarray*}  
\lefteqn{-{d\over ds}\Lambda^{(m)}_n(\bar\theta_{1,n};\,\bar x_{1,n}-s) 
= {1\over\mu}\sum_{i=1}^n\theta_i\,  
 \Xi^{(m)}_n(\bar\theta_{1,n};\,\bar x_{1,n}-s)}\\ 
&&-e^{-\theta_1(x_1-s)/\mu}\left(
 {d\over ds} \Lambda^{(m)}_{n-1}(\bar\theta_{2,n};\bar x_{2,n}-s)
+ {1\over\mu}\sum_{i=2}^n\theta_i\, \Xi^{(m)}_{n-1}(\bar\theta_{2,n};\bar x_{2,n}-s)\right)\\
&&-{1\over m}{d\over ds} R_1^{(m)}(\bar x_{1,n}-s)
 -{1\over m}{1\over\mu}\sum_{i=1}^n \,\theta_i R^{(m)}(\bar x_{1,n}-s).
\end{eqnarray*}
Since 
\begin{eqnarray*}
R_1^{(m)}(\bar x_{1,n}-s)
={\theta_1\over \mu}
  \int_s^{x_1} e^{-\theta_1(x_1-x)/\mu}
   \,\widetilde\Lambda^{(m)}_{n-1}(\bar\theta_{2,n};\bar x_{2,n}-x)\,
      \widetilde\Lambda^{(m)}_1(\eta_n;\,x-s) 
\,dx
\end{eqnarray*}
we have the identity
\begin{eqnarray*}
&&-{d\over ds} R_1^{(m)}(\bar x_{1,n}-s)\\
&&\quad ={\theta_1\over \mu}
  \int_s^{x_1} e^{-\theta_1(x_1-x)/\mu}
   \,\widetilde\Lambda^{(m)}_{n-1}(\bar\theta_{2,n};\bar x_{2,n}-x)\,
         d_x\widetilde\Lambda^{(m)}_1(\eta_n;\,x-s)
\\
&&\quad ={1\over\mu}\sum_{i=1}^n \theta_i\, R^{(m)}(\bar x_{1,n}-s).
\end{eqnarray*}
Thus,
\begin{eqnarray*}  
\lefteqn{{d\over ds}\Lambda^{(m)}_n(\bar\theta_{1,n};\,\bar x_{1,n}-s) 
+\frac{\eta_n}{\mu}\,  
 \Xi^{(m)}_n(\bar\theta_{1,n};\,\bar x_{1,n}-s)}\\ 
&&=e^{-\theta_1(x_1-s)/\mu}\left(
 {d\over ds} \Lambda^{(m)}_{n-1}(\bar\theta_{2,n};\bar x_{2,n}-s)
+ {1\over\mu}\sum_{i=2}^n\theta_i\, \Xi^{(m)}_{n-1}(\bar\theta_{2,n};\bar x_{2,n}-s)\right).
\end{eqnarray*}
The statement of the lemma now follows by induction.
\hfill $\Box$

\begin{lemma}\label{lemmaXi} For each $n\ge 1$, the limit functions 
\[
\Xi_n(\bar\theta_{1,n};\,\bar x_{1,n})=\lim_{m\to\infty}
\Xi_n^{(m)}(\bar\theta_{1,n};\,\bar x_{1,n})
\] 
exist and are given by 
\begin{eqnarray}\nonumber
\lefteqn{\Xi_n(\bar\theta_{1,n};\,\bar x_{1,n})
=\sum_{j=1}^n \exp\Big\{-
   \sum_{i=1}^{j-1}(\theta_i+\dots+\theta_n)(x_i-x_{i-1})/\mu\Big\}}\\
 &&\quad\times{1\over\mu}\int_{x_{j-1}}^{x_j}(\theta_j+\dots+\theta_n)
    e^{-(\theta_j+\dots+\theta_n)(u-x_{j-1})/\mu} \beta^{-1}u^{-\beta}\,du. 
\label{xiformula}
\end{eqnarray}
For $n\ge 2$ they solve the recursive system
\[
\Xi_n(\bar\theta_{1,n};\,\bar x_{1,n}) 
- \Xi_1(\eta_n;\,x_1)  
= e^{-\theta_1x_1/\mu}\Big(
\Xi_{n-1}(\bar\theta_{2,n};\bar x_{2,n})        
- \Xi_1(\eta_n-\theta_1;\,x_1)  
\Big). 
\]
\end{lemma} 
\textbf{Proof}. 
As $m\to\infty$, by Lemmas \ref{lemmaincreasing} and \ref{lemmaTminusTtilde},
\begin{equation}\label{xi_1}
\Xi_1^{(m)}(\theta;\,x)={\mu\over\theta}
{d\over dx}mE(e^{\theta(T_{ax}-ax/\mu)/a}-1)
\to{1\over\beta\mu}\int_0^x \theta e^{-\theta u/\mu}u^{-\beta}\,du
=\Xi_1(\theta;\,x).
\end{equation}
Using (\ref{derivativeonedim}), for arbitrary $\alpha$, 
\begin{eqnarray*}\nonumber
\lefteqn{\int_0^{x_1} e^{-\theta_1(x_1-x)/\mu}\,
            \Xi_1^{(m)}(\alpha;\,dx)}\\
&&=\Xi_1^{(m)}(\alpha;\,x_1)
   - {\theta_1\over\mu} \int_0^{x_1} e^{-\theta_1(x_1-x)/\mu}\,
            \Xi_1^{(m)}(\alpha;\,x)\,dx\\
&&=\Xi_1^{(m)}(\alpha;\,x_1)
   - {\theta_1\over\alpha}
\int_0^{x_1} e^{-\theta_1(x_1-x)/\mu}\,
     \Lambda_1^{(m)}(\alpha;\,dx) \label{Xi1bound}\\
&&=\Xi_1^{(m)}(\alpha;\,x_1)
   - {\theta_1\over\alpha}\Lambda_1^{(m)}(\alpha;\,x_1)
+{\theta_1^2\over\mu\alpha}\int_0^{x_1} e^{-\theta_1(x_1-x)/\mu}\,
     \Lambda_1^{(m)}(\alpha;\,x)\,dx.
\end{eqnarray*}
Here, (\ref{xi_1}) shows that $\Xi_1^{(m)}$ converges to $\Xi_1$ and
Lemma \ref{lemmamargscaling} shows that $\Lambda_1^{(m)}$ converges to
a limit function $\Lambda_1$ which satisfies ${d\over dx}
\Lambda_1(\alpha,x)=\alpha\,\Xi_1(\alpha,x)/\mu$. Since $0\le
\Lambda_1^{(m)}(\alpha;\,x)\le \Lambda_1^{(m)}(\alpha;\,x_1)$ for
$0\le x\le x_1$ by Lemma \ref{lemmaincreasing}, we can find a
dominating function for $\Lambda_1^{(m)}(\alpha;\,x)$ on $[0,x_1]$ and
conclude that the last integral term also converges. By reverting the
partial integrations this yields with $\alpha=\eta_n$ 
\begin{equation}
\int_0^{x_1} e^{-\theta_1(x_1-x)/\mu}\,
            \Xi_1^{(m)}(\eta_n;\,dx) 
\to 
\int_0^{x_1} e^{-\theta_1(x_1-x)/\mu}\,
            \Xi_1(\eta_n;\,dx). 
\end{equation}

Let us now consider the remainder terms $R^{(m)}$ in Lemma
\ref{lemmaXi_m}.  Because of (\ref{derivativeonedim}) we may rewrite using 
\begin{eqnarray*}
\widetilde\Lambda_1^{(m)}(\alpha;\,dx)
&=& \Lambda_1^{(m)}(\alpha;\,dx)+\Xi_1^{(m)}(\alpha;\,dx)\\
&=& \Lambda_1^{(m)}(\alpha;\,dx)
 +e^{-\alpha x/\mu}\,d_x(e^{\alpha x/\mu}\Xi_1^{(m)}(\alpha;\,x))
  -{\alpha\over\mu}\, \Xi_1^{(m)}(\alpha;\,x)\,dx\\
&=& e^{-\alpha x/\mu}\,d_x(e^{\alpha x/\mu}\Xi_1^{(m)}(\alpha;\,x)),
\end{eqnarray*}
and obtain, again with $\alpha=\sum_{i=1}^n\theta_i$,
\[
R^{(m)}= {\theta_1\over \alpha}
e^{-\theta_1 x_1/\mu} \int_0^{x_1} e^{-(\alpha-\theta_1)x/\mu} 
   \,\widetilde\Lambda^{(m)}_{n-1}(\bar\theta_{2,n};\bar x_{2,n}-x)\,
  d_x\Big(e^{\alpha x/\mu}\Xi^{(m)}_1(\alpha;x)\Big).
\]
The above integration is carried out with respect to the function
\[
F^{(m)}(x)=e^{\alpha x/\mu}\,\Xi^{(m)}_1(\alpha;\,x)
=m E(e^{\alpha \widetilde T_{ax}/a}-e^{\alpha T_{ax}/a}).
\]
It was observed in the proof of Lemma \ref{lemmaincreasing} that $T$ can
be viewed as a shift of $\widetilde T$ with the first point of
increase in $X_0$. In particular, $\widetilde T_x-T_x$ equals $\widetilde
T_{X_0}$ almost surely on the set $\{X_0<x\}$. Hence the increment
\begin{eqnarray*}
\lefteqn{E(e^{\alpha \widetilde T_{x+h}}-e^{\alpha T_{x+h}})
 -E(e^{\alpha \widetilde T_x}-e^{\alpha T_x})}\\
&&=E((e^{\alpha T_{x+h}}-e^{\alpha T_x})(e^{\alpha\widetilde
  T_{X_0}}-1),\, X_0<x)\\
&&\quad + E(e^{\alpha T_{x+h}}(e^{\alpha \widetilde T_{X_0}} 
-1)
-(e^{\alpha \widetilde T_x}-1),\,x<X_0<x+h)\\
&&\quad  +E(e^{\alpha \widetilde T_{x+h}}
 -e^{\alpha \widetilde T_x},\,x+h<X_0)
\end{eqnarray*}
is positive for $\alpha>0$ and negative for $\alpha<0$.  Thus,
$F^{(m)}$ is a monotone measure with limit $e^{\alpha
x/\mu}\,\Xi_1(\alpha;\,x)$. Using the variation measure
$|F^{(m)}|$ we obtain a constant $C_n(x_1)$ uniform in $m$ such that
\begin{eqnarray*}
|R^{(m)}| &\le&  \Big|{\theta_1\over \alpha}\Big|
  \sup_{0\le x\le x_1}\Big| e^{-(\theta_1 x_1+(\alpha-\theta_1)x)/\mu}
   \,\widetilde\Lambda^{(m)}_{n-1}(\bar\theta_{2,n}; \bar x_{2,n}-x)\Big|\,
   |F^{(m)}|(x_1)\\
&\le & C_n(x_1)\, \sup_{0\le x\le x_1}|
\widetilde\Lambda^{(m)}_{n-1}(\bar\theta_{2,n}; \bar x_{2,n}-x)|
\end{eqnarray*}
The same arguments apply to the remainder terms $R_1^{(m)}$ in
(\ref{Lambdaeqn}) and $R_2^{(m)}$ in (\ref{tildeLambdaeqn}). Hence,
taking $C_n(x_1)$ sufficiently large,
\[
R_i^{(m)}\le  C_n(x_1)\, \sup_{0\le x\le x_1}|
\widetilde\Lambda^{(m)}_{n-1}(\bar\theta_{2,n}; \bar x_{2,n}-x)|,\quad
i=1,2.
\]

We are now prepared to carry out an induction on $n$ in Lemma
\ref{lemmaXi_m} and equation (\ref{tildeLambdaeqn}). Assume that 
$\Xi^{(m)}_{n-1}(\bar\theta_{2,n}; \bar x_{2,n})$ converges to 
$\Xi_{n-1}(\bar\theta_{2,n}; \bar x_{2,n})$  
and $\widetilde\Lambda^{(m)}_{n-1}(\bar\theta_{2,n}; \bar x_{2,n})$
is such that 
\[
\sup_{0\le x\le x_1}|
\widetilde\Lambda^{(m)}_{n-1}(\bar\theta_{2,n}; \bar
x_{2,n}-x)|<C_{n-1}(\bar x_{1,n}).
\]
Then, applying (\ref{xi_1}), (\ref{Xi1bound}) and the induction
hypothesis to Lemma \ref{lemmaXi_m} it follows that all limit functions
$\Xi_n$ exist and satisfy
\begin{eqnarray*}\nonumber
\lefteqn{\quad \Xi_n(\bar\theta_{1,n};\,\bar x_{1,n}) 
- \Xi_1(\eta_n;\,x_1)  
  \;=\; e^{-\theta_1x_1/\mu}\,
\Xi_{n-1}(\bar\theta_{2,n};\bar x_{2,n})}\\
&& \quad + \Big(\frac{\theta_1}{\eta_n}-1\Big)
    \int_0^{x_1} e^{-\theta_1(x_1-x)/\mu}\,
            \Xi_1(\eta_n;\,dx) 
\\
&&\!\!\!\!=\;  e^{-\theta_1x_1/\mu}\Big(
\Xi_{n-1}(\bar\theta_{2,n};\bar x_{2,n})        
- {1\over\mu}(\eta_n-\theta_1) 
\int_0^{x_1} e^{-x(\eta_n-\theta_1)/\mu}   
      \beta^{-1}x^{-\beta}\,dx\Big),
\end{eqnarray*}
which is the desired relation. Moreover, observing
that the convergence of $\widetilde I_4^{(m)}$ is a byproduct of the
proof of (\ref{Xi1bound}), it follows from (\ref{tildeLambdaeqn}) that 
\[
\sup_{0\le u\le x_0}|
\widetilde\Lambda^{(m)}_n(\bar\theta_{1,n}; \bar
x_{1,n}-u)|<C_n(\bar x_{0,n}),\quad 0<x_0<x_1<\dots < x_n.
\]

To verify the explicit form
(\ref{xiformula}) of the solution, assume that the
claim is correct for index $n-1$. Then
\begin{eqnarray*}
\lefteqn{\Xi_{n-1}(\bar\theta_{2,n};\bar x_{2,n}) 
= {1\over\mu}\int_0^{x_2}(\theta_2+\dots+\theta_n)  
    e^{-(\theta_2+\dots+\theta_n)u/\mu}
    \beta^{-1}u^{-\beta}\,du}\\
&&+\sum_{j=2}^{n-1} \exp\Big\{-(\theta_2 +\dots+\theta_n)x_2/\mu
 -\sum_{i=2}^{j-1}(\theta_{i+1}+\dots+\theta_n)(x_{i+1}-x_{i})/\mu\Big\}\\
&&\quad\times{1\over\mu}\int_{x_{j}}^{x_{j+1}}(\theta_{j+1}+\dots+\theta_n)
    e^{-(\theta_{j+1}+\dots+\theta_n)(u-x_{j})/\mu}
    \beta^{-1}u^{-\beta}\,du.
\end{eqnarray*}
This implies
\begin{eqnarray*}
\lefteqn{\Xi_{n-1}(\bar\theta_{2,n};\bar x_{2,n}) 
-\Xi_1\Big(\sum_{i=2}^n\theta_i\,;\,x_1\Big)}\\
&&= {1\over\mu}\int_{x_1}^{x_2}(\theta_2+\dots+\theta_n)
    e^{-(\theta_2+\dots+\theta_n)u/\mu}
    \beta^{-1}u^{-\beta}\,du\\
&&+\sum_{j=3}^{n} \exp\Big\{-(\theta_2 +\dots+\theta_n)x_2/\mu
 -\sum_{i=3}^{j-1}(\theta_i+\dots+\theta_n)(x_i-x_{i-1})/\mu\Big\}\\
&&\quad\times{1\over\mu}\int_{x_{j-1}}^{x_{j}}(\theta_{j}+\dots+\theta_n)
    e^{-(\theta_{j}+\dots+\theta_n)(u-x_{j-1})/\mu}
    \beta^{-1}u^{-\beta}\,du.
\end{eqnarray*}
Hence 
\begin{eqnarray*}
\lefteqn{e^{-\theta_1x_1/\mu}\Big(\Xi_{n-1}(\bar\theta_{2,n};\bar
  x_{2,n}) -\Xi_1\Big(\sum_{i=2}^n\theta_i\,;\,x_1\Big)\Big)}\\ 
&&=e^{-(\theta_1+\dots+\theta_n)x_1/\mu} {1\over\mu}\int_{x_1}^{x_2}
(\theta_2+\dots+\theta_n)
  e^{-(\theta_2+\dots+\theta_n)(u-x_1)/\mu} \beta^{-1}u^{-\beta}\,du\\
&&\quad+\sum_{j=3}^{n} \exp\Big\{
  -\sum_{i=1}^{j-1}(\theta_i+\dots+\theta_n)(x_i-x_{i-1})/\mu\Big\}\\
&&\qquad\times{1\over\mu}\int_{x_{j-1}}^{x_{j}}(\theta_{j}+\dots+\theta_n)
  e^{-(\theta_{j}+\dots+\theta_n)(u-x_{j-1})/\mu}
  \beta^{-1}u^{-\beta}\,du\\
&&=\sum_{j=2}^{n} \exp\Big\{
  -\sum_{i=1}^{j-1}(\theta_i+\dots+\theta_n)(x_i-x_{i-1})/\mu\Big\}\\
&&\quad\times{1\over\mu}\int_{x_{j-1}}^{x_{j}}(\theta_{j}+\dots+\theta_n)
  e^{-(\theta_{j}+\dots+\theta_n)(u-x_{j-1})/\mu}
  \beta^{-1}u^{-\beta}\,du.
\end{eqnarray*}
 \hfill $\Box$

The remaining proofs of the convergence of multivariate distributions
in Theorem 2 are organized in three consecutive lemmas, leading up to
the identification of the cumulant generating function (\ref{mgfndim})
in Theorem 2. 

\begin{lemma}\label{lemmalambdaeqn}
The limit functions $\Lambda_n=\lim_{m\to\infty}\Lambda_n^{(m)}$,
$n\ge 1$, exist and we have 
\begin{equation}\label{lambda_1}
\Lambda_1(\theta;x)
 ={\theta^2\over\beta\mu^2}
  \int_0^x\int_0^u e^{-\theta v/\mu}v^{-\beta}\,dvdu
\end{equation}
and for $n\ge 2$, recursively
\begin{eqnarray*}
\Lambda_n(\bar\theta_{1,n};\,\bar x_{1,n}) 
&=&\Lambda_{n-1}(\bar\theta_{2,n};\bar x_{2,n}-x_1)
       +\Lambda_1(\eta_n;\,x_1) 
\\ 
&& +\Big(\frac{\theta_1}{\eta_n}-1\Big)
    \int_0^{x_1} e^{-\theta_1(x_1-x)/\mu} 
            \Lambda_1(\eta_n;\,dx)
\\
&& -\Big(\frac{\theta_1}{\eta_n-\theta_1}+1\Big)
   \int_0^{x_1} e^{-\theta_1(x_1-x)/\mu}
   \,d_x\Lambda_{n-1}(\bar\theta_{2,n};\bar x_{2,n}-x). 
\end{eqnarray*}
\end{lemma}
\textbf{Proof}. For $n=1$ this follows from Lemma \ref{lemmamargscaling}
and for $n\ge 2$ from (\ref{Lambdaeqn}) and a further partial
integration of the term $I^{(m)}_2$, which gives
\begin{eqnarray*}
\lefteqn{\Lambda_n^{(m)}(\bar\theta_{1,n};\,\bar x_{1,n}) 
=\Lambda^{(m)}_{n-1}(\bar\theta_{2,n};\bar x_{2,n}-x_1)
       +\Lambda^{(m)}_1(\eta_n;\,x_1) 
}\\ 
&&\quad +\Big(\frac{\theta_1}{\eta_n}-1\Big)
    \int_0^{x_1} e^{-\theta_1(x_1-x)/\mu} 
            \Lambda^{(m)}_1(\eta_n;\,dx)
\\
&&\quad -\int_0^{x_1} e^{-\theta_1(x_1-x)/\mu}
   \,d_x\Lambda^{(m)}_{n-1}(\bar\theta_{2,n};\bar x_{2,n}-x)\\
&&\quad +{\theta_1\over\mu}\int_0^{x_1} e^{-\theta_1(x_1-x)/\mu}\,
   \Xi^{(m)}_{n-1}(\bar\theta_{2,n};\bar x_{2,n}-x)\,dx
         +{1\over m}R_1^{(m)}\\
&&=\Lambda^{(m)}_{n-1}(\bar\theta_{2,n};\bar x_{2,n}-x_1)
       +\Lambda^{(m)}_1(\eta_n;\,x_1) 
\\ 
&&\quad +\Big(\frac{\theta_1}{\eta_n}-1\Big)
    \int_0^{x_1} e^{-\theta_1(x_1-x)/\mu} 
            \Lambda^{(m)}_1(\eta_n;\,dx) 
\\
&&\quad -\Big(\frac{\theta_1}{\eta_n-\theta_1}+1\Big)
\int_0^{x_1} e^{-\theta_1(x_1-x)/\mu}
   \,d_x\Lambda^{(m)}_{n-1}(\bar\theta_{2,n};\bar x_{2,n}-x)
         +{1\over m}R_1^{(m)},
\end{eqnarray*}
where we apply Lemma \ref{lemmaderivative_m} for the last equality.  The
arguments which justify that we are allowed to exchange the order of
integration and taking limits in $m$, as well as controlling the
remainder terms, are parallel to those in the proof of Lemma
\ref{lemmaXi}, again based on Lemma \ref{lemmaderivative_m}. \hfill $\Box$

\medskip

In view of (\ref{defmgf}) we conclude from Lemma \ref{lemmalambdaeqn} the
convergence of the finite-dimensional distributions in Theorem \ref{thmICR}. 

\begin{lemma} \label{lemmalambdafinal}
The finite-dimensional distributions of the sequence of random
processes studied in Theorem \ref{thmICR} (with $c=1$) converge to
those of a limit process $Y_\beta$, such that the collection of
logarithmic moment generating functions
\[
\Lambda_n(\bar\theta_{1,n};\,\bar x_{1,n}) =\log E
 \exp\Big\{\sum_{i=1}^n\theta_i Y_\beta(x_i)\Big\},\quad n\ge 1 
\]
is the unique solution to the closed system of linear integral
equations
\begin{eqnarray*}
\lefteqn{\Lambda_n(\bar\theta_{1,n};\,\bar x_{1,n}) 
=\Lambda_{n-1}(\bar\theta_{2,n};\bar x_{2,n}-x_1)
       +\Lambda_1\Big(\sum_{i=1}^n\theta_i;\,x_1\Big)}\\ 
&& +\Big({\theta_1\over\sum_{i=1}^n\theta_i}-1\Big)
    \int_0^{x_1} e^{-\theta_1(x_1-x)/\mu} 
            \Lambda_1\Big(\sum_{i=1}^n\theta_i;\,dx\Big)\\
&& -\Big({\theta_1\over\sum_{i=2}^n\theta_i}+1\Big)
\int_0^{x_1} e^{-\theta_1(x_1-x)/\mu}
   \,d_x  \Lambda_{n-1}(\bar\theta_{2,n};\bar x_{2,n}-x),\quad n\ge 2,
\end{eqnarray*}
with $\Lambda_1$ as in (\ref{lambda_1}). 
\end{lemma}

\begin{lemma} The cumulant function for the increments of $Y_\beta$,
\[
\Gamma_n(\bar\theta_{1,n};\,\bar x_{1,n}) =\log E
 \exp\Big\{\sum_{i=1}^n\theta_i(Y_\beta(x_i)-Y_\beta(x_{i-1}))\Big\}
\]
has the explicit form given in (\ref{mgfndim}).
\end{lemma}
\textbf{Proof}. We have
\[
\Gamma_n(\bar\theta_{1,n};\,\bar x_{1,n})  
=\Lambda_n((\theta_1-\theta_2,\dots,\theta_{n-1}-\theta_n,\theta_n),\bar
x_{1,n})
\]
so by Lemma \ref{lemmalambdafinal} 
\begin{eqnarray*}
\lefteqn{\Gamma_n(\bar\theta_{1,n};\,\bar x_{1,n}) 
=\Gamma_{n-1}(\bar\theta_{2,n};\bar x_{2,n}-x_1)
       +\Lambda_1(\theta_1;\,x_1)}\\ 
&& -{\theta_2\over\theta_1}
     \int_0^{x_1} e^{-(\theta_1-\theta_2)(x_1-x)/\mu} 
            \Lambda_1(\theta_1;\,dx)\\
&& -{\theta_1\over\theta_2}
  \int_0^{x_1} e^{-(\theta_1-\theta_2)(x_1-x)/\mu}
   \,d\Gamma_{n-1}(\bar\theta_{2,n};\bar x_{2,n}-x),\quad n\ge 2.
\end{eqnarray*}
It may now be checked that the functions in (\ref{mgfndim}) solve the
above system of equations. For details, see Gaigalas, Kaj
\cite{gaigalaskaj}, Section 6.3. \hfill $\Box$

\section{Remaining proofs}

\subsection{Limiting distribution under FBM scaling}

In this section we discuss briefly the convergence of the finite-dimensional
distributions in Theorem 1. Recall that for standard fractional
Brownian motion $B_H$,
\[
\log E \exp\Big\{\sum_{i=1}^n\theta_i \sigma_\beta B_H(x_i)\Big\}=
  {1\over 2}\sigma^2_\beta\sum_{i=1}^n\sum_{j=1}^n\theta_i\theta_j\, 
      {\rm Cov}(B_H(x_i),B_H(x_j))
\]
with
\[
{\rm Cov}(B_H(x),B_H(y))={1\over 2}(x^{2-\beta}+y^{2-\beta}-(x-y)^{2-\beta}).
\]
In the scaling regime defined by (\ref{FCR}) and (\ref{FCR_b}) we have
\[
{a\over b}= \sqrt{{a^\beta\mu\over mL(a)}}\to 0,\qquad
{am\over b}= \sqrt{{a^\beta m\mu\over L(a)}}\to \infty.
\]
By analyzing in this case the recursive equations for 
$\widetilde\Lambda^{(m)}_n(\bar\theta_{1,n};\bar x_{1,r})$ it follows
that
\begin{equation}\label{FCRlambda}
\widetilde\Lambda_n(\bar\theta_{1,n};\bar x_{1,n})=\lim_{m\to\infty}
\widetilde\Lambda^{(m)}_n(\bar\theta_{1,n};\bar x_{1,r})
={\delta^2_\beta\over\mu}\sum_ {i=1}^n\theta_i x_i^{1-\beta}, \quad
\delta^2_\beta={1\over \beta(1-\beta)}. 
\end{equation}
Moreover, the limit functions 
$\Lambda_n(\bar\theta_{1,n};\bar x_{1,r})=
\lim_{m\to\infty}\Lambda^{(m)}_n(\bar\theta_{1,n};\bar x_{1,r})$ 
satisfy, in analogy to the result of Lemma \ref{lemmalambdafinal},
\[
\Lambda_1(\theta;\,x)=\lim_{m\to\infty}\Lambda_1^{(m)}(\theta;\,x)
  = {1\over 2}\sigma_\beta^2\mu^{-2} \theta^2 x^{2-\beta}
\]
and for $n\ge 2$, 
\begin{eqnarray*}
\Lambda_n(\bar\theta_{1,n};\,\bar x_{1,n}) 
&=&\Lambda_{n-1}(\bar\theta_{2,n};\,\bar x_{2,n})
 +{\theta_1\over \sum_{i=1}^n\theta_i}
       \Lambda_1\Big(\sum_{i=1}^n\theta_i;\,x_1\Big)\\
&&+{\theta_1\over \mu}\int_0^{x_1}
   \widetilde\Lambda_{n-1}(\bar\theta_{2,n};\bar x_{2,n}-x)\,dx.
\end{eqnarray*}
Thus, using (\ref{FCRlambda}), 
\[
\Lambda_n(\bar\theta_{1,n};\,\bar x_{1,n}) 
=\Lambda_{n-1}(\bar\theta_{2,n};\,\bar x_{2,n})
+{\sigma_\beta^2\over 2\mu^2} \sum_{j=1}^n \theta_1\theta_j
  {1\over 2}\Big[x_1^{2-\beta}+x_j^{2-\beta}-(x_j-x_1)^{2-\beta}\Big].
\]
Hence
\begin{eqnarray*} 
\Lambda_n(\bar\theta_{1,n};\,\bar x_{1,n}) 
={\sigma^2_\beta\over 2\mu^2}\sum_{i=1}^n
  \sum_{j=1}^n \theta_i\theta_j {1\over 2}(x_i^{2-\beta}+x_j^{2-\beta}
    -(x_i-x_j)^{2-\beta}).
\end{eqnarray*}

\subsection{Proof of tightness in C}

To complete the proofs of our results we establish tightness of
the sequences 
\[
Y^{(m)}(x)={1\over b_m}\sum_{i=1}^m(T^{(i)}_{a_mx}-{a_mx\over\mu})
\]
studied in Theorems 1 and 2, by applying a standard moment criterion. 
Since $Y^{(m)}(x)$ has stationary increments, to prove that
$\{Y^{(m)}\}$ is tight in $C$ it is enough to find $\gamma>1$, an
integer $m_0$ and a constant $K$ such that for fixed $T$,
\begin{equation}\label{criteriontight}
{\rm Var}(Y^{(m)}(x))={m\over b_m^2} {\rm Var}(T_{a_m x}) \le K x^\gamma
\end{equation}
for $0<x<T$ and $m\ge m_0$ (Billingsley (1968), Thm.\ 12.3). 

By Lemma \ref{lemmaTtildemean} iii), the variance of $T_x$ is a
non-decreasing function in $x$.  Hence we may apply Karamata's
Tauberian theorem (Bingham {\it et al.} \cite{binghametal} Theorem
1.7.1) to show that ${\rm Var}(T_x)$ is regularly varying in infinity
with index $2-\beta$. Indeed, recalling the previously used
notation $I(u)=\mu u -\Phi(u)$, the asymptotic property
(\ref{I_karamata}) implies
\[  
\int_0^\infty u e^{-ux} {\rm Var}(T_x)\,dx
 = {2\over(\mu u)^2}\sum_{n=1}^\infty \Big({I(u)\over\mu u}\Big)^n
\sim {2\Gamma(1-\beta) L(1/u)\over\beta\mu^2 u^{2-\beta}}
,\quad u\to 0, 
\]
hence
\begin{equation}\label{varT_regvar}
{\rm Var}(T_x)\sim 
{2\Gamma(1-\beta) x^{2-\beta}L(x)\over\Gamma(3-\beta)\beta\mu^3}
={\sigma^2_\beta \over\mu^3} L(x)x^{2-\beta} ,\quad x\to\infty.
\end{equation}
The next step is to apply the Potter bounds for regularly varying
functions (Bingham {\it et al.}, Ch 1.5) to obtain for any
$\epsilon>0$ an $a_0$, such that
\[
{{\rm Var}(T_{ax})\over {\rm Var}(T_a)}\le (1+\epsilon)
\max(x^{2-\beta+\epsilon},x^{2-\beta-\epsilon}),\quad a\ge a_0,\; ax\ge a_0.
\]
Hence for $m\ge m_0$ so large that $a\ge a_0$, $ax\ge a_0$,
\[
m{\rm Var}(T_{ax})/b^2\le (1+\epsilon)m{\rm Var}(T_a)b^{-2}
\max(x^{2-\beta+\epsilon},x^{2-\beta-\epsilon}) 
\]
But in either case of the FBM scaling (\ref{FCR}), (\ref{FCR_b}) or
the intermediate scaling (\ref{ICR}), the asymptotic relation
(\ref{varT_regvar}) yields
\[
m{\rm Var}(T_{a})/b^2\to
    \sigma^2_\beta/\mu^2,\quad m\to\infty,
\]
and so, eventually choosing a larger $m_1\ge m_0$,
\[
m{\rm Var}(T_{ax})/b^2\le (1+\epsilon)(\sigma^2_\beta+\epsilon)
\max(x^{2-\beta+\epsilon},x^{2-\beta-\epsilon}), \quad m\ge m_1.
\]
With $\epsilon<1-\beta$ this yields (\ref{criteriontight}) for $ax\ge a_0$. 

It remains to prove (\ref{criteriontight}) for $a\ge a_0$ and
$ax<a_0$.  By Lemma \ref{lemmaTtildemean} iv),
\[
m{\rm Var}(T_{ax})/b^2 
\le {2e\over\mu} {m\over b^2}\int_0^{ax}\Phi(1/y)^{-1}\,dy
\le {2e\over\mu} {m\over b^2} \frac{ax}{\Phi(1/ax)}
\]
By (\ref{lowerindex}) we can find a constant $C_1$ and $q>\beta$, such
that $\phi(1/ax)\le C_1 (ax)^q$.  As in the proof of Lemma
\ref{lemmaTminusTtilde}, we may take $a^{-\epsilon}\le L(a)$.
In Theorem 1, $ma^{2-\beta}L(a)/\mu b^2=1$. In Theorem 2,
$ma^{2-\beta}L(a)/\mu b^2\to c^\beta$.  Thus,
\[
m{\rm Var}(T_{ax})/b^2 
\le C_2 {ma^{2-\eta}L(a)\over b^2}
\frac{(ax)^{1+q}}{a^{2-\beta-\epsilon}}
\le C_3 \, a_0^\beta x^{1+q-\beta}\frac{1}{a^{1-q-\epsilon}} 
\]
Since $\beta<q<\sigma\le 1$, we may take $\epsilon<1-q$ to obtain 
(\ref{criteriontight}) for $\gamma=1+q-\beta>1$. 

\hfill $\Box$


\end{document}